# NONPARAMETRIC ESTIMATION OF SCALAR DIFFUSIONS BASED ON LOW FREQUENCY DATA[1]

By Emmanuel Gobet, Marc Hoffmann and Markus Reiß

*Ecole Polytechnique, Université Paris 7 and Humboldt-Universität zu Berlin*

We study the problem of estimating the coefficients of a diffusion $(X_t, t \geq 0)$; the estimation is based on discrete data $X_{n\Delta}, n = 0, 1, \ldots, N$. The sampling frequency $\Delta^{-1}$ is constant, and asymptotics are taken as the number $N$ of observations tends to infinity. We prove that the problem of estimating both the diffusion coefficient (the volatility) and the drift in a nonparametric setting is ill-posed: the minimax rates of convergence for Sobolev constraints and squared-error loss coincide with that of a, respectively, first- and second-order linear inverse problem. To ensure ergodicity and limit technical difficulties we restrict ourselves to scalar diffusions living on a compact interval with reflecting boundary conditions.

Our approach is based on the spectral analysis of the associated Markov semigroup. A rate-optimal estimation of the coefficients is obtained via the nonparametric estimation of an eigenvalue–eigenfunction pair of the transition operator of the discrete time Markov chain $(X_{n\Delta}, n = 0, 1, \ldots, N)$ in a suitable Sobolev norm, together with an estimation of its invariant density.

## 1. Introduction.

1.1. *Overview.* Since Feller's celebrated classification, stationary scalar diffusions have served as a representative model for homogeneous Markov processes in continuous time. Historically, diffusion processes were probably first seen as approximation models for discrete Markov chains, up to an appropriate rescaling in time and space. More recently, the development of financial mathematics has argued in favor of genuinely continuous time

Received July 2002; revised June 2003.
[1]Supported in part by the European network Dynstoch and the DFG-Sonderforschungsbereich 373.

models, with simple dynamics governed by a local mean (drift) $b(\cdot)$ and local variance (diffusion coefficient, or volatility) $\sigma(\cdot)$ on the state space $S = \mathbb{R}$ or $S \subset \mathbb{R}$ with appropriate boundary condition. The dynamics are usually described by an Itô-type stochastic differential equation in the interior of $S$, which in the time-homogeneous case reads like

$$dX_t = b(X_t)\,dt + \sigma(X_t)\,dW_t, \qquad t \geq 0,$$

where the driving process $(W_t, t \geq 0)$ is standard Brownian motion. The growing importance of diffusion models progressively raised among the community of statisticians a vast research program, from both quantitative and theoretical angles. We outline the main achievements of this program in Section 1.2.

In the late 1970s a statistician was able to characterize qualitatively the properties of a parametric ergodic diffusion model based on the continuous observation of a sample path

$$X^T := (X_t, 0 \leq t \leq T)$$

of the trajectory, as $T \to \infty$, that is, as the time length of the experiment grows to infinity, a necessary assumption to assure the growing of information thanks to the recurrence of the sample path. The 1980s explored various discretization schemes of the continuous time model: the data $X^T$ could progressively be replaced by the more realistic observation

$$X^{(N,\Delta_N)} := (X_{n\Delta_N}, n = 0, 1, \ldots, N),$$

with asymptotics taken as $N \to \infty$. The discretization techniques used at that time required the high frequency sampling assumption $\Delta_N \to 0$ whereas $N\Delta_N \to \infty$ in order to guarantee the closeness of $X^{(N,\Delta_N)}$ and $X^T$, with $T = N\Delta_N$. Soon, a similar nonparametric program was achieved for both continuous time and high frequency data.

By the early to mid-1990s, the frontier remained the "fixed $\Delta$ case," that is, the case of low frequency data. This is the topic of the present paper. First, one must understand the importance and flexibility gained by being able to relax the assumption that the sampling time $\Delta$ between two data points is "small": indeed, one can hardly deny that, in practice, it may well happen that sampling with arbitrarily small $\Delta$ is simply not feasible. Put differently, the asymptotic statistical theory is a mathematical construct to assess the quality of an estimator based on discrete observations and it must be decided which asymptotics are adequate for the data at hand. Second, the statistical nature of the problem drastically changes when passing from high to low frequency sampling: the approximation properties of the sample path $X^{N\Delta_N}$ by $X^{(N,\Delta_N)}$ are not valid anymore; the observation $(X_0, X_\Delta, \ldots, X_{N\Delta})$ becomes a genuine Markov chain, and inference about the underlying coefficients of the diffusion process must be sought via the identification of the law of the observation $X^{(N,\Delta_N)}$. In the time-homogeneous case the mathematical properties of the random vector $X^{(N,\Delta_N)}$ are embodied in the transition operator

$$P_\Delta f(x) := \mathbb{E}[f(X_\Delta)|X_0 = x],$$

defined on appropriate test functions $f$. Under suitable assumptions, the operator $P_\Delta$ is associated with a Feller semigroup $(P_t, t \geq 0)$ with a densely defined infinitesimal generator $L$ on the space of continuous functions given by

$$Lf(x) = L_{\sigma,b}f(x) := \frac{\sigma^2(x)}{2}f''(x) + b(x)f'(x).$$

The second-order term $\sigma(\cdot)$ is the diffusion coefficient, and the first-order term $b(\cdot)$ is the drift coefficient. Postulating the existence of an invariant density $\mu(\cdot) = \mu_{\sigma,b}(\cdot)$, the operator $L$ is unbounded, but self-adjoint negative on $L^2(\mu) := \{f | \int |f|^2 \mu < \infty\}$, and the functional calculus gives the correspondence

$$(1.1) \qquad P_\Delta = \exp(\Delta L)$$

in the operator sense. Therefore, a consistent statistical program can be presumed to start from the observed Markov chain $X^{(N,\Delta)}$, estimate its transition operator $P_\Delta$ and infer about the pair $(b(\cdot), \sigma(\cdot))$, via the correspondence (1.1), in other words via the spectral properties of the operator $P_\Delta$. Expressed in a diagram, we obtain the following line:

$$(1.2) \qquad \text{data} = X^{(N,\Delta)} \xrightarrow{(E)} P_\Delta \xrightarrow{(I)} L \longleftrightarrow (b(\cdot), \sigma(\cdot)) = \text{parameter}.$$

The efficiency of a given statistical estimation procedure will be measured by the proficiency in combining the estimation part $(E)$ and the identification part $(I)$ of the model.

The works of Hansen, Scheinkman and Touzi (1998) and Chen, Hansen and Scheinkman (1997) paved the way: they formulated a precise and thorough program, proposing and discussing several methods for identifying scalar diffusions via their spectral properties. Simultaneously, the Danish school, given on impulse by the works of Kessler and Sørensen (1999), systematically studied the parametric efficiency of spectral methods in the fixed $\Delta$ setting described above. By constructing estimating functions based on eigenfunctions of the operator $L$, they could construct $\sqrt{N}$-consistent estimators and obtained precise asymptotic properties.

However, a quantitative study of nonparametric estimation in the fixed $\Delta$ context remained out of reach for some time, for both technical and conceptual reasons. The purpose of the present paper is to fill in this gap, by trying to understand and explain why the nonparametric case significantly differs from its parametric analogue, as well as from the high frequency data framework in nonparametrics.

We are going to establish minimax rates of convergence over various smoothness classes, characterizing upper and lower bounds for estimating $b(\cdot)$ and $\sigma(\cdot)$ based on the obervation of $X_0, X_\Delta, \ldots, X_{N\Delta}$, with asymptotics taken as $N \to \infty$. The minimax rate of convergence is an index of both accuracy of estimation and complexity of the model. We will show that in the nonparametric case the complexity of the problems of estimating $b(\cdot)$ and $\sigma(\cdot)$ is related to ill-posed inverse problems. Although we mainly focus on the theoretical aspects of the statistical model, the estimators we propose are based on feasible nonparametric smoothing methods: they can be implemented in practice, allowing for adaptivity and finite sample optimization. Some simulation results were performed by Reiß (2003).

The estimation problem is exactly formulated in Section 2, where also the main theoretical results are stated. The spectral estimation method we adopt is explained in Section 3, which includes a discussion of related problems and possible extensions. The proofs of the upper bound for our estimator and its optimality in a minimax sense are given in Sections 4 and 5, respectively. Results of rather technical nature are deferred to Section 6.

1.2. *Statistical estimation for diffusions: an outlook.* We give a brief and selective summary of the evolution of the area over the last two decades. The nonparametric identification of diffusion processes from continuous data was probably first addressed in the reference paper of Banon (1978). More precise estimation results can be listed as follows:

1.2.1. *Continuous or high frequency data: the parametric case.* Estimation of a finite-dimensional parameter $\theta$ from $X^T = (X_t, 0 \leq t \leq T)$ with asymptotics as $T \to \infty$ when $X$ is a diffusion of the form

$$dX_t = b_\theta(X_t)\, dt + \sigma(X_t)\, dW_t \tag{1.3}$$

is classical [Brown and Hewitt (1975) and Kutoyants (1975)]. Here $(W_t, t \geq 0)$ is a standard Wiener process. The diffusion coefficient is perfectly identified from the data by means of the quadratic variation of $X$. By assuming the process $X$ to be ergodic (positively recurrent), a sufficiently regular parametrization $\theta \mapsto b_\theta(\cdot)$ implies the local asymptotic normality (LAN) property for the underlying statistical model, therefore ensuring the $\sqrt{T}$-consistency and efficiency of the ML-estimator [see Liptser and Shiryaev (2001)].

In the case of discrete data $X_{n\Delta_N}, n = 0, 1, \ldots, N$, with high frequency sampling $\Delta_N^{-1} \to \infty$, but long range observation $N\Delta_N \to \infty$ as $N \to \infty$, various discretization schemes and estimating procedures had been proposed [Yoshida (1992) and Kessler (1997)] until Gobet (2002) eventually proved the LAN property for ergodic diffusions of the form

$$dX_t = b_{\theta_1}(X_t)\,dt + \sigma_{\theta_2}(X_t)\,dW_t \tag{1.4}$$

in a general setting, by means of the Malliavin calculus: under suitable regularity conditions, the finite-dimensional parameter $\theta_1$ in the drift term can be estimated with optimal rate $\sqrt{N\Delta_N}$, whereas the finite-dimensional parameter $\theta_2$ in the diffusion coefficient is estimated with the optimal rate $\sqrt{N}$.

1.2.2. *Continuous or high frequency data*: *the nonparametric case.* A similar program was progressively obtained in nonparametrics: If the drift function $b(\cdot)$ is globally unknown in the model given by (1.3), but belongs to a Sobolev ball $S(s, L)$ (of smoothness order $s > 0$ and radius $L$) over a given compact interval $\mathcal{I}$, a certain kernel estimator $\hat{b}_T(\cdot)$ achieves the following upper bound in $L^2(\mathcal{I})$ and in a root-mean-squared sense:

$$\sup_{b \in S(s,L)} \mathbb{E}[\|\hat{b}_T - b\|^2_{L^2(\mathcal{I})}]^{1/2} \lesssim T^{-s/(2s+1)}.$$

This already indicates a formal analogy with the model of nonparametric regression or "signal + white noise" where the same rate holds. (Here and in the sequel, the symbol $\lesssim$ means "up to constants," possibly depending on the parameters of the problem, but that are continuous in their arguments.) See Kutoyants (1984) for precise mathematical results.

Similar extensions to the discrete case with high frequency data sampling for the model driven by (1.4) were given in Hoffmann (1999), where the rates $(N\Delta_N)^{-s/(2s+1)}$ for the drift function $b(\cdot)$ and $N^{-s/(2s+1)}$ for the diffusion coefficient $\sigma(\cdot)$ have been obtained and proved to be optimal. See also the pioneering paper of Pham (1981). Methods of this kind have been successfully applied to financial data [Aït-Sahalia (1996), Stanton (1997), Chapman and Pearson (2000) and Fan and Zhang (2003)]. In particular, it is investigated whether the usual parametric model assumptions are compatible with the data, and the use of nonparametric methods is advocated.

1.2.3. *From high to low frequency data.* As soon as the sampling frequency $\Delta_N^{-1} = \Delta^{-1}$ is not large anymore, the problem of estimating a parameter in the drift or diffusion coefficient becomes significantly more difficult: the trajectory properties that can be recovered from the data when $\Delta_N$ is small are lost. In particular, there is no evident approximating scheme that can efficiently compute or mimic the continuous ML-estimator in parametric estimation.

Likewise, the usual nonparametric kernel estimators, based on differencing, do not provide consistent estimation of the drift $b(\cdot)$ or the diffusion coefficient $\sigma(\cdot)$. As a concrete example, consider the standard Nadaraya–Watson estimator $\hat{b}(x)$ of the drift $b(x)$ in the point $x \in \mathbb{R}$:

$$\hat{b}(x) := \frac{(N\Delta)^{-1} \sum_{n=0}^{N-1} K_h(x - X_{n\Delta})(X_{(n+1)\Delta} - X_{n\Delta})}{N^{-1} \sum_{n=0}^{N-1} K_h(x - X_{n\Delta})}$$

with a kernel function $K(\cdot)$ and $K_h(x) := h^{-1}K(h^{-1}x)$ for $h > 0$. If we let $N \to \infty$ and $h \to 0$, then by the ratio ergodic theorem and by kernel properties we obtain almost surely the limit

$$\mathbb{E}[\Delta^{-1}(X_\Delta - X_0)|X_0 = x] = \Delta^{-1} \int_0^\Delta P_t b(x)\, dt.$$

Hence, this estimator is not consistent. It merely yields a highly blurred version of $b(x)$, which of course tends to $b(x)$ in the high frequency limit $\Delta \to 0$. Note that the transition operators $P_t$ involved depend on the unknown functions $b(\cdot)$ and $\sigma(\cdot)$ as a whole. The situation for estimators of $\sigma(\cdot)$ based on the approximation of the quadratic variation is even worse, because the drift $b(\cdot)$ enters directly into the limit expression.

1.2.4. *Spectral methods for parametric estimation.* Kessler and Sørensen (1999) suggested the use of eigenvalues $\kappa_\theta$ and eigenvectors $\varphi_\theta(\cdot)$ of the parametrized infinitesimal generator

$$L_\theta f(x) = \frac{\sigma_\theta^2(x)}{2} f''(x) + b_\theta(x) f'(x),$$

that is, such that $L_\theta \varphi_\theta(x) = \kappa_\theta \varphi_\theta(x)$. Indeed, since the pair $(\kappa_\theta, \varphi_\theta)$ also satisfies

$$P_\Delta \varphi_\theta(X_{n\Delta}) = \mathbb{E}[\varphi_\theta(X_{(n+1)\Delta})|X_{n\Delta}] = \exp(\kappa_\theta \Delta)\varphi_\theta(X_{n\Delta}),$$

whenever it is easy to compute, the knowledge of a pair $(\kappa_\theta, \varphi_\theta)$ can be translated into a set of conditional moment conditions to be used in estimating functions. With their method, Kessler and Sørensen can construct $\sqrt{N}$-consistent estimators that are nearly efficient. See also the paper of Hansen, Scheinkman and Touzi (1998) that we already mentioned.

In a sense, in this idea also lies the essence of our method. However, the strategy of Kessler and Sørensen is not easily extendable to nonparametrics: there is no straightforward way to pass from a finite-dimensional parametrization of the generator $L_\theta$ with explicit eigenpairs $(\kappa_\theta, \varphi_\theta)$ to a full nonparametric space with satisfactory approximation properties. Besides, there would be no evident route to treat the variance of such speculative nonparametric estimators either, because the behavior of the parametric Fisher information matrix for a growing number of parameters is too complex to be easily controlled. We will see in Section 3 how to pass over these objections by estimating directly an eigenpair nonparametrically.

1.2.5. *Prospectives.*  A quick summary yields Table 1 for optimal rates of convergence.

Table 1 can be interpreted as follows: the difficulty of the estimation problem is increasing from top to bottom and from left to right. A blank line separates the continuous–high-frequency (HF) data domain from the low-frequency (LF) data domain. The breach for LF data opened by Kessler and Sørensen as well as by Hansen, Scheinkman and Touzi shows that $\sqrt{N}$-consistent estimators usually exist in the parametric case. The remaining case are the rates of convergence for LF data in the nonparametric case $u_N$ for the drift $b(\cdot)$ and $v_N$ for the diffusion coefficient $\sigma(\cdot)$, for which we are aiming.

## 2. Main results.

2.1. *A diffusion model with boundary reflections.*  We shall restrict ourselves to reflecting diffusions on a one-dimensional interval to avoid highly nontrivial technical issues; see the discussion in Section 3.3.

Choosing for convenience the interval $[0,1]$, we suppose the following.

ASSUMPTION 2.1.  The function $b:[0,1] \to \mathbb{R}$ is measurable and bounded, the function $\sigma:[0,1] \to (0,\infty)$ is continuous and positive and the function $\nu:[0,1] \to \mathbb{R}$ satisfies $\nu(0)=1$, $\nu(1)=-1$.

We consider the stochastic differential equation

$$dX_t = b(X_t)\,dt + \sigma(X_t)\,dW_t + \nu(X_t)\,dL_t(X),$$
(2.1)
$$X_0 = x_0 \quad \text{and} \quad X_t \in [0,1] \qquad \forall t \geq 0.$$

The process $(W_t, t \geq 0)$ is a standard Brownian motion and $(L_t(X), t \geq 0)$ is a nonanticipative continuous nondecreasing process that increases only when $X_t \in \{0,1\}$. The boundedness of $b(\cdot)$ and the ellipticity of $\sigma(\cdot)$ ensure the existence of a weak solution; see for instance Stroock and Varadhan (1971). Note that the process $L(X)$ is part of the solution and is given by a difference of local times of $X$ at the boundary points of $[0,1]$.

TABLE 1

|  | Parametric | | Nonparametric | |
|---|---|---|---|---|
|  | $b$ | $\sigma$ | $b$ | $\sigma$ |
| Continuous | $T^{-1/2}$ | known | $T^{-s/(2s+1)}$ | known |
| HF data | $(N\Delta_N)^{-1/2}$ | $N^{-1/2}$ | $(N\Delta_N)^{-s/(2s+1)}$ | $N^{-s/(2s+1)}$ |
| LF data | $N^{-1/2}$ | $N^{-1/2}$ | $u_N$ | $v_N$ |

Due to the compactness of $[0,1]$ and the reflecting boundary conditions, the Markov process $X$ has a spectral gap, which implies geometric ergodicity; compare with Lemmas 6.1 and 6.2. In particular, a unique invariant measure $\mu$ exists and the one-dimensional distributions of $X_t$ converge exponentially fast to $\mu$ as $t \to \infty$ so that the assumption of stationarity can be made without loss of generality for asymptotic results.

We denote by $\mathbb{P}_{\sigma,b}$ the law of the associated stationary diffusion on the canonical space $\Omega = \mathcal{C}(\mathbb{R}_+, [0,1])$ of continuous functions over the positive axis with values in $[0,1]$, equipped with the topology of uniform convergence and endowed with its Borel $\sigma$-field $\mathcal{F}$. We denote by $\mathbb{E}_{\sigma,b}$ the corresponding expectation operator. Given $N \geq 1$ and $\Delta > 0$, we observe the canonical process $(X_t, t \geq 0)$ at equidistant times $n\Delta$ for $n = 0, 1, \ldots, N$. Let $\mathcal{F}_N$ denote the $\sigma$-field generated by $\{X_{n\Delta} | n = 0, \ldots, N\}$.

DEFINITION 2.2. An estimator of the pair $(\sigma(\cdot), b(\cdot))$ is an $\mathcal{F}_N$-measurable function on $\Omega$ with values in $L^2([0,1]) \times L^2([0,1])$.

To assess the $L^2$-risk in a minimax framework, we introduce the nonparametric set $\Theta_s$, which consists of pairs of functions of regularity $s$ and $s-1$, respectively.

DEFINITION 2.3. For $s > 1$ and given constants $C \geq c > 0$, we consider the class $\Theta_s := \Theta(s, C, c)$ defined by

$$\left\{ (\sigma, b) \in H^s([0,1]) \times H^{s-1}([0,1]) \big| \|\sigma\|_{H^s} \leq C,\ \|b\|_{H^{s-1}} \leq C,\ \inf_x \sigma(x) \geq c \right\},$$

where $H^s$ denotes the $L^2$-Sobolev space of order $s$.

Note that all $(\sigma(\cdot), b(\cdot)) \in \Theta_s$ satisfy Assumption 2.1.

2.2. *Minimax rates of convergence.* We are now in position to state the main theorems. By (3.12) and (3.13) in the next section we define estimators $\hat{\sigma}^2$ and $\hat{b}$ using a spectral estimation method based on the observation $(X_0, X_\Delta, \ldots, X_{N\Delta})$. These estimators, which of course depend on the number $N$ of observations, satisfy the following uniform asymptotic upper bounds.

THEOREM 2.4. *For all $s > 1$, $C \geq c > 0$ and $0 < a < b < 1$ we have*

$$\sup_{(\sigma,b) \in \Theta_s} \mathbb{E}_{\sigma,b}[\|\hat{\sigma}^2 - \sigma^2\|_{L^2([a,b])}^2]^{1/2} \lesssim N^{-s/(2s+3)},$$

$$\sup_{(\sigma,b) \in \Theta_s} \mathbb{E}_{\sigma,b}[\|\hat{b} - b\|_{L^2([a,b])}^2]^{1/2} \lesssim N^{-(s-1)/(2s+3)}.$$

Recall that $A \lesssim B$ means that $A$ can be bounded by a constant multiple of $B$, where the constant depends continuously on other parameters involved. Similarly, $A \gtrsim B$ is equivalent to $B \lesssim A$ and $A \sim B$ holds if both relations $A \lesssim B$ and $A \gtrsim B$ are true.

As the following lower bounds prove, the rates of convergence of our estimators are optimal in a minimax sense over $\Theta_s$.

THEOREM 2.5. *Let $E_N$ denote the set of all estimators according to Definition* 2.2. *Then for all $0 \leq a < b \leq 1$ and $s > 1$ the following hold:*

$$(2.2) \qquad \inf_{\hat{\sigma}^2 \in E_N} \sup_{(\sigma,b) \in \Theta_s} \mathbb{E}_{\sigma,b}[\|\hat{\sigma}^2 - \sigma^2\|^2_{L^2([a,b])}]^{1/2} \gtrsim N^{-s/(2s+3)},$$

$$(2.3) \qquad \inf_{\hat{b} \in E_N} \sup_{(\sigma,b) \in \Theta_s} \mathbb{E}_{\sigma,b}[\|\hat{b} - b\|^2_{L^2([a,b])}]^{1/2} \gtrsim N^{-(s-1)/(2s+3)}.$$

If we set $s_1 = s - 1$ and $s_2 = s$, then the drift $b(\cdot) \in \Theta_s$ with regularity $s_1$ can be estimated with the minimax rate of convergence $u_N = N^{-s_1/(2s_1+5)}$, whereas the diffusion coefficient $\sigma(\cdot) \in \Theta_s$ has regularity $s_2$ and the corresponding minimax rate of convergence is $v_N = N^{-s_2/(2s_2+3)}$. Hence, Table 1 in Section 1.2.5 can be filled with two rather unexpected rates of convergence $u_N$ and $v_N$. In Section 3.3 the rates are explained in the terminology of ill-posed problems and reasons are given why the tight connection between the regularity assumptions on $b(\cdot)$ and $\sigma(\cdot)$ is needed.

## 3. Spectral estimation method.

3.1. *The basic idea.* We shall base our estimator of the diffusion coefficient $\sigma(\cdot)$ and of the drift coefficient $b(\cdot)$ on spectral methods for passing from the transition operator $P_\Delta$, which is approximately known to us, to the infinitesimal generator $L$, which more explicitly encodes the functions $\sigma(\cdot)$ and $b(\cdot)$. In the sequel, we shall rely on classical results for scalar diffusions; for example, consult Bass [(1998), Chapter 4]. We use the specific form of the invariant density

$$(3.1) \qquad \mu(x) = 2C_0 \sigma^{-2}(x) \exp\left(\int_0^x 2b(y)\,\sigma^{-2}(y)\,dy\right)$$

and the function $S(\cdot) = 1/s'(\cdot)$, derived from the scale function $s(\cdot)$,

$$(3.2) \qquad S(x) = \tfrac{1}{2}\sigma^2(x)\mu^{-1}(x) = C_0 \exp\left(-\int_0^x 2b(y)\sigma^{-2}(y)\,dy\right),$$

with the normalizing constant $C_0 > 0$ depending on $\sigma(\cdot)$ and $b(\cdot)$. The action of the generator in divergence form is given by

$$(3.3) \quad Lf(x) = L_{\sigma,b}f(x) = \frac{1}{2}\sigma^2(x)f''(x) + b(x)f'(x) = \frac{1}{\mu(x)}(S(x)f'(x))',$$

where the domain of this unbounded operator on $L^2(\mu)$ is given by the subspace of the $L^2$-Sobolev space $H^2$ with Neumann boundary conditions

$$\mathrm{dom}(L) = \{f \in H^2([0,1]) | f'(0) = f'(1) = 0\}.$$

The generator $L$ is a self-adjoint elliptic operator on $L^2(\mu)$ with compact resolvent so that it has nonpositive point spectrum only. If $\nu_1$ denotes the largest negative eigenvalue of $L$ with eigenfunction $u_1$, then due to the reflecting boundary of $[0,1]$ the Neumann boundary conditions $u_1'(0) = u_1'(1) = 0$ hold and thus

$$(3.4) \quad Lu_1 = \mu^{-1}(Su_1')' = \nu_1 u_1 \quad \Longrightarrow \quad S(x)u_1'(x) = \nu_1 \int_0^x u_1(y)\mu(y)\,dy.$$

From (3.2) we can derive an explicit expression for the diffusion coefficient:

$$(3.5) \qquad \sigma^2(x) = \frac{2\nu_1 \int_0^x u_1(y)\mu(y)\,dy}{u_1'(x)\mu(x)}.$$

The corresponding expression for the drift coefficient is

$$(3.6) \qquad b(x) = \nu_1 \frac{u_1(x)u_1'(x)\mu(x) - u_1''(x)\int_0^x u_1(y)\mu(y)\,dy}{u_1'(x)^2\mu(x)}.$$

Hence, if we knew the invariant measure $\mu$, the eigenvalue $\nu_1$ and the eigenfunction $u_1$ (including its first two derivatives), we could exactly determine the drift and diffusion coefficient. Of course, these identities are valid for any eigenfunction $u_k$ with eigenvalue $\nu_k$, but for better numerical stability we shall use only the largest nondegenerate eigenvalue $\nu_1$. Moreover, it is known that only the eigenfunction $u_1$ does not have a vanishing derivative in the interior of the interval (cf. Proposition 6.5) so that by this choice indeterminacy at interior points is avoided.

Using semigroup theory [Engel and Nagel (2000), Theorem IV.3.7] we know that $u_1$ is also an eigenfunction of $P_\Delta$ with eigenvalue $\kappa_1 = e^{\Delta \nu_1}$. Our procedure consists of determining estimators $\hat{\mu}$ of $\mu$ and $\hat{P}_\Delta$ of $P_\Delta$, to calculate the corresponding eigenpair $(\hat{\kappa}_1, \hat{u}_1)$ and to use (3.5) and (3.6) to build a plug-in estimator of $\sigma(\cdot)$ and $b(\cdot)$.

3.2. *Construction of the estimators.* We use projection methods, taking advantage of approximating properties of abstract operators by finite-dimensional matrices, for which the spectrum is easy to calculate numerically. A similar approach was already suggested by Chen, Hansen and Scheinkman (1997). More specifically, we make use of wavelets on the interval $[0,1]$. For the construction of wavelet bases and their properties we refer to Cohen (2000).

DEFINITION 3.1. Let $(\psi_\lambda)$ with multiindices $\lambda = (j,k)$ be a compactly supported $L^2$-orthonormal wavelet basis of $L^2([0,1])$. The approximation spaces $(V_J)$ are defined as $L^2$-closed linear spans of the wavelets up to the frequency level $J$,

$$V_J := \overline{\mathrm{span}}\{\psi_\lambda | |\lambda| \leq J\} \qquad \text{where } |(j,k)| := j.$$

The $L^2$-orthogonal projection onto $V_J$ is called $\pi_J$; the $L^2(\mu)$-orthogonal projection onto $V_J$ is called $\pi_J^\mu$.

In the sequel we shall regularly use the Jackson and Bernstein inequalities with respect to the $L^2$-Sobolev spaces $H^s([0,1])$ of regularity $s$:

$$\|(\mathrm{Id} - \pi_J)f\|_{H^t} \lesssim 2^{-J(s-t)} \|f\|_{H^s}, \qquad 0 \leq t \leq s,$$

$$\forall v_J \in V_J, \qquad \|v_J\|_{H^s} \lesssim 2^{J(s-t)} \|v_J\|_{H^t}, \qquad 0 \leq t \leq s.$$

The canonical projection estimate of $\mu$ based on $(X_{n\Delta})_{0 \leq n \leq N}$ is given by

$$(3.7) \qquad \hat{\mu} := \sum_{|\lambda| \leq J} \hat{\mu}_\lambda \psi_\lambda \qquad \text{with } \hat{\mu}_\lambda := \frac{1}{N+1} \sum_{n=0}^{N} \psi_\lambda(X_{n\Delta}).$$

By the ergodicity of $X$ it follows that $\hat{\mu}_\lambda$ is a consistent estimate of $\langle \mu, \psi_\lambda \rangle$ for $N \to \infty$. To estimate the action of the transition operator on the wavelet basis $(\mathbf{P}_\Delta^J)_{\lambda,\lambda'} := \langle P_\Delta \psi_\lambda, \psi_{\lambda'} \rangle_\mu$, we introduce the symmetrized matrix estimator $\hat{\mathbf{P}}_\Delta$ with entries

$$(3.8)\ (\hat{\mathbf{P}}_\Delta)_{\lambda,\lambda'} := \frac{1}{2N} \sum_{n=1}^{N} (\psi_\lambda(X_{(n-1)\Delta})\psi_{\lambda'}(X_{n\Delta}) + \psi_{\lambda'}(X_{(n-1)\Delta})\psi_\lambda(X_{n\Delta})).$$

This yields an approximation of $\langle P_\Delta \psi_{\lambda'}, \psi_\lambda \rangle_\mu$, that is, of the action of the transition operator on $V_J$ with respect to the unknown scalar product $\langle \cdot, \cdot \rangle_\mu$ in $L^2(\mu)$. We therefore introduce a third statistic $\hat{\mathbf{G}}$, which approximates the $\dim(V_J) \times \dim(V_J)$-dimensional Gram matrix $\mathbf{G}$ with entries $\mathbf{G}_{\lambda,\lambda'} = \langle \psi_\lambda, \psi_{\lambda'} \rangle_\mu$, and which is given by

$$(3.9) \quad \hat{\mathbf{G}}_{\lambda,\lambda'} := \frac{1}{N}\left(\frac{1}{2}\psi_\lambda(X_0)\psi_{\lambda'}(X_0) + \frac{1}{2}\psi_\lambda(X_{N\Delta})\psi_{\lambda'}(X_{N\Delta}) + \sum_{n=1}^{N-1} \psi_\lambda(X_{n\Delta})\psi_{\lambda'}(X_{n\Delta})\right).$$

The particular treatment of the boundary terms will be explained later. If we put $\boldsymbol{\Sigma} = (w_n \psi_\lambda(X_n))_{|\lambda| \leq J, n \leq N}$ with $w_0 = w_N = \frac{1}{2}$ and $w_n = 1$ otherwise, we have $\hat{\mathbf{G}} = N^{-1} \boldsymbol{\Sigma}\boldsymbol{\Sigma}^T$, $\boldsymbol{\Sigma}^T$ being the transpose of $\boldsymbol{\Sigma}$. Our construction can thus be regarded as a least squares type estimator, as in a usual regression setting; see the argument developed in Section 3.3.1.

We combine the last two estimators in order to determine estimates for the eigenvalue $\kappa_1$ and the eigenfunction $u_1$ of $P_\Delta$. As will be made precise in Proposition 4.5, the operators $P_\Delta$ and $\pi_J^\mu P_\Delta$ are close for large values of $J$. Note that all eigenvectors of $\pi_J^\mu P_\Delta$ lie in $V_J$, the range of $\pi_J^\mu P_\Delta$. The eigenfunction $u_1^J$ corresponding to the second largest eigenvalue $\kappa_1^J$ of $\pi_J^\mu P_\Delta$ is characterized by

$$\langle P_\Delta u_1^J, \psi_\lambda \rangle_\mu = \kappa_1^J \langle u_1^J, \psi_\lambda \rangle_\mu \qquad \forall |\lambda| \leq J. \tag{3.10}$$

We pass to vector-matrix notation and use from now on bold letters to define for a function $v \in V_J$ the corresponding coefficient column vector $\mathbf{v} = (\langle v, \psi_\lambda \rangle)_{|\lambda| \leq J}$. Observe carefully the different $L^2$-scalar products used; here they are with respect to the Lebesgue measure. Thus, we can rewrite (3.10) as

$$\mathbf{P}_\Delta^J \mathbf{u}_1^J = \kappa_1^J \mathbf{G} \mathbf{u}_1^J. \tag{3.11}$$

As $\mathbf{v}^T \mathbf{G} \mathbf{v} = \langle v, v \rangle_\mu > 0$ holds for $v \in V_J \setminus \{0\}$, the matrix $\mathbf{G}$ is invertible and $(\kappa_1^J, \mathbf{u}_1^J)$ is an eigenpair of $\mathbf{G}^{-1} \mathbf{P}_\Delta^J$. This matrix is self-adjoint with respect to the scalar product induced by $\mathbf{G}$:

$$\langle \mathbf{G}^{-1} \mathbf{P}_\Delta^J \mathbf{v}, \mathbf{w} \rangle_\mathbf{G} := (\mathbf{G}^{-1} \mathbf{P}_\Delta^J \mathbf{v})^T \mathbf{G} \mathbf{w}$$
$$= \mathbf{v}^T \mathbf{P}_\Delta^J \mathbf{w} = \langle \mathbf{v}, \mathbf{G}^{-1} \mathbf{P}_\Delta^J \mathbf{w} \rangle_\mathbf{G}.$$

Similarly, $\mathbf{v}^T \hat{\mathbf{G}} \mathbf{v} = N^{-1} (\mathbf{\Sigma}^T \mathbf{v})^T \mathbf{\Sigma}^T \mathbf{v} \geq 0$ holds and the matrix $\hat{\mathbf{G}}$ can be shown to be even strictly positive definite with high probability (see Lemma 4.12). In this case, we similarly infer that $\hat{\mathbf{G}}^{-1} \hat{\mathbf{P}}$ is self-adjoint with respect to the $\hat{\mathbf{G}}$-scalar product. The Cauchy–Schwarz inequality and the inequality between geometric and arithmetic mean yield the estimate

$$\langle \hat{\mathbf{G}}^{-1} \hat{\mathbf{P}}_\Delta \mathbf{v}, \mathbf{v} \rangle_{\hat{\mathbf{G}}} = \frac{1}{N} \sum_{n=1}^N v(X_{(n-1)\Delta}) v(X_{n\Delta})$$
$$\leq \frac{1}{N} \left( \sum_{n=0}^{N-1} v(X_{n\Delta})^2 \right)^{1/2} \left( \sum_{n=1}^N v(X_{n\Delta})^2 \right)^{1/2}$$
$$\leq \frac{1}{N} \left( \frac{1}{2} v(X_0)^2 + \frac{1}{2} v(X_{N\Delta})^2 + \sum_{n=1}^{N-1} v(X_{n\Delta})^2 \right)$$
$$= \langle \mathbf{v}, \mathbf{v} \rangle_{\hat{\mathbf{G}}}.$$

We infer that all eigenvalues of $\hat{\mathbf{G}}^{-1} \hat{\mathbf{P}}_\Delta$ are real and not larger than 1. Hence, the second largest eigenvalue $\hat{\kappa}_1$ of $\hat{\mathbf{G}}^{-1} \hat{\mathbf{P}}_\Delta$ is well defined, which is why we downweighted the boundary terms of $\hat{\mathbf{G}}$. The eigenvalue $\hat{\kappa}_1$ of $\hat{\mathbf{G}}^{-1}\hat{\mathbf{P}}_\Delta$ and its corresponding eigenvector $\hat{\mathbf{u}}_1$ yield estimators of $\kappa_1^J$ and $u_1^J$.

Plugging the estimator $\hat{\mu}$ as well as $\hat{\kappa}_1$ and $\hat{u}_1$ into (3.5) and (3.6), we obtain our estimators of $\sigma^2(\cdot)$ and $b(\cdot)$:

$$(3.12) \quad \hat{\sigma}^2(x) := \frac{2\Delta^{-1}\log(\hat{\kappa}_1)\int_0^x \hat{u}_1(y)\hat{\mu}(y)\,dy}{\hat{u}_1'(x)\hat{\mu}(x)},$$

$$(3.13) \quad \hat{b}(x) := \Delta^{-1}\log(\hat{\kappa}_1)\frac{\hat{u}_1(x)\hat{u}_1'(x)\hat{\mu}(x) - \hat{u}_1''(x)\int_0^x \hat{u}_1(y)\hat{\mu}(y)\,dy}{\hat{u}_1'(x)^2\hat{\mu}(x)}.$$

To avoid indeterminacy, the denominators of the estimators are forced to remain above a certain minimal level, which depends on the subinterval $[a,b] \subset [0,1]$ for which the loss function is taken. See (4.5) for the exact formulation in the case of $\hat{\sigma}^2(\cdot)$ and proceed analogously for $\hat{b}(\cdot)$.

### 3.3. Discussion.

3.3.1. *Least squares approach.* The estimator matrix $\hat{\mathbf{G}}^{-1}\hat{\mathbf{P}}_\Delta$ is built as in the least squares approach for projection methods in classical regression. To estimate $P_\Delta \psi_{\lambda_0}(x) = \mathbb{E}_{\sigma,b}[\psi_{\lambda_0}(X_\Delta)|X_0 = x]$, the least squares method consists of minimizing

$$(3.14) \quad \sum_{n=1}^N \left| \psi_{\lambda_0}(X_{n\Delta}) - \sum_{|\lambda|\leq J} \alpha_\lambda^0 \psi_\lambda(X_{(n-1)\Delta}) \right|^2 \longrightarrow \min!$$

over all real coefficients $(\alpha_\lambda^0)$, leading to the normal equations

$$\sum_{n=1}^N \left( \sum_{|\lambda|\leq J} \alpha_\lambda^0 \psi_\lambda(X_{(n-1)\Delta}) \right) \psi_{\lambda'}(X_{(n-1)\Delta}) = \sum_{n=1}^N \psi_{\lambda'}(X_{(n-1)\Delta})\psi_{\lambda_0}(X_{n\Delta})$$

for all $|\lambda'| \leq J$. Up to the special treatment of the boundary terms, we thus obtain the vector $(\alpha_\lambda^0)$ as the column with index $\lambda_0$ in $\hat{\mathbf{G}}^{-1}\hat{\mathbf{P}}_\Delta$.

3.3.2. *Other than wavelet methods.* For our projection estimates to work, we merely need approximation spaces satisfying the Jackson and Bernstein inequalities. Hence, other finite element bases could serve as well.

The invariant density and the transition density could also be estimated using kernel methods, but the numerical calculation of the eigenpair $(\hat{\kappa}_1, \hat{u}_1)$ would then involve an additional discretization step.

3.3.3. *Diffusions over the real line.* Using a wavelet basis of $L^2(\mathbb{R})$, it is still possible to estimate $\mu$ and $P_\Delta$ over the real line; in particular the eigenvalue characterization (3.10) extends to this case. Hansen, Scheinkman and Touzi (1998) derive the same formulae as (3.5) and (3.6) under ergodicity and boundary conditions so that a plug-in approach is feasible. However, a theoretical study seems to require much more demanding theoretical tools. If the uniform separation of the spectral value $\nu_1$ and a polynomial growth bound for the eigenfunction $u_1(\cdot)$ are ensured, we expect that the same minimax results hold with respect to an $L^2(\mu)$-loss function, where the invariant density $\mu(\cdot)$ is of course parameter-dependent. However, all spectral approximation results have to be reconsidered with extra care, in particular because the $L^2(\mu)$-norms are in general not equivalent for different parameters.

3.3.4. *Adaptation to unknown smoothness.* The knowledge of the smoothness $s$ that is needed for the construction of our estimators is not realistic in practice. An adaptive estimation of the eigenpair $(u_1(\cdot), \kappa_1)$ and $\mu(\cdot)$ that yields adaptive estimators for $(\sigma(\cdot), b(\cdot))$ could be obtained by the following modifications: First, the adaptive estimation of $\mu(\cdot)$ in a classical mixing framework is fairly well known [e.g., Tribouley and Viennet (1998)]. Second, taking advantage of the multiresolution structure provided by wavelets, the adaptive estimation of $P_\Delta$ could be obtained by introducing an appropriate thresholding in the estimated matrices on a large approximation space.

3.3.5. *Interpretation as an ill-posed problem.* One can make the link with ill-posed inverse problems by saying that estimation of $\mu(\cdot)$ is well-posed (i.e., with achievable rate $N^{-s/(2s+1)}$), but for $S(\cdot)$ we need an estimate of the derivative $u_1'(\cdot)$ yielding an ill-posedness degree of 1 ($N^{-s/(2s+3)}$). Observe that the regularity conditions $\sigma \in H^s$ and $b \in H^{s-1}$ are translated into $\mu \in H^s$, $S \in H^s$. The transformation of $(\mu, S)$ to $\sigma^2(\cdot) = 2S(\cdot)/\mu(\cdot)$ is stable [$L^2$-continuous for $S(\cdot) \geq s_0 > 0$], whereas in $b(\cdot) = S'(\cdot)/\mu(\cdot)$ another ill-posed operation (differentiation) occurs with degree 1.

A brief stepwise explanation reads as follows. Step 1, the natural parametrization $(\mu, P_\Delta)$ is well-posed (for $P_\Delta$ in the strong operator norm sense). Step 2, the calculation of the spectral pair $(\kappa_1, u_1)$ is well-posed. Step 3, the differentiation of $u_1$ that determines $S$ has an ill-posedness of degree 1. Step 4, the calculation of $\sigma^2$ from $(\mu, S)$ is well-posed. Step 5, the calculation of $b$ from $(\mu, S)$ is ill- posed of degree 1.

3.3.6. *Regularity restrictions on $b(\cdot)$ and $\sigma(\cdot)$.* It is noteworthy that in the continuous time or high frequency observation case, the parameter $b(\cdot)$ does not influence the asymptotic behavior of the estimator of $\sigma(\cdot)$ and vice versa. The estimation problems are separated. In our low frequency regime we had to suppose tight regularity connections between $\sigma(\cdot)$ and $b(\cdot)$. This stems from the fact that for the underlying Markov chain $X^{(N,\Delta)}$ the parameters $\mu(\cdot)$ and $S(\cdot)$ are more natural and the regularity of these functions depends on the regularity both of $b(\cdot)$ and of $\sigma(\cdot)$.

At a different level, in nonparametric regression, different smoothness constraints are needed between the mean and the variance function. Recommended references are Müller and Stadtmüller (1987) and Fan and Yao (1998).

Finally, although we ask for the tight connection $s_1 = s_2 - 1$ for the regularity $s_1$ of the drift $b(\cdot)$ and $s_2$ of the diffusion coefficient $\sigma(\cdot)$, our results readily carry over to the milder constraint $s_1 \geq s_2 - 1$.

3.3.7. *Estimation when one parameter is known.* If $\sigma(\cdot)$ is known, an estimate $\hat{\mu}$ of the invariant density yields an estimate of $b(\cdot)$, since

$$b(x) = \frac{(\sigma^2(x)\mu(x))'}{2\mu(x)}, \qquad x \in [0,1].$$

Estimation of $\mu \in H^s$, $s > 1$, in $H^1$-norm can be achieved with rate $N^{-(s-1)/(2s+1)}$ and this rate is thus also valid for estimating $b(\cdot)$ in $L^2$-norm. Given the drift coefficient $b(\cdot)$, we find

$$\sigma^2(x) = 2\frac{\int_0^x b(y)\mu(y)\,dy + C}{\mu(x)}, \qquad x \in [0,1],$$

where $C$ is a suitable constant. If we knew $\sigma^2(0)$, we would obtain the rate $N^{-s/(2s+1)}$ for $\mu \in H^s$.

Using a preliminary nonparametric estimate $\hat{\sigma}_C^2$ depending on the parameter $C$ and then fitting a parametric model for $C$, we are likely to find the same rate. In any case, the assumption of knowing one parameter seems rather artificial and no further investigations have been performed.

3.3.8. *Estimation at the boundary.* Our plug-in estimators can only be defined on the open subinterval $(0,1)$. Estimation at the boundary points leads to a risk increase due to $S^{-1}(0) = \nu_1 u_1(0)\mu(0)/u_1''(0)$ by de l'Hôspital's rule applied to (3.4). Thus, estimating $\sigma(0)$ and $b(0)$ involves taking the second and third derivative, respectively, when using plug-in estimators. A pointwise lower bound result—along the lines of the $L^2$-lower bound proof—shows that this deterioration cannot be avoided.

## 4. Proof of the upper bound.

4.1. *Convergence of $\hat{\mu}$.* First, we recall the proof for the risk bound in estimating the invariant measure:

PROPOSITION 4.1. *With the choice $2^J \sim N^{1/(2s+1)}$ the following uniform risk estimate holds for $\hat{\mu}$ based on $N$ observations:*

$$\sup_{(\sigma,b)\in\Theta_s} \mathbb{E}_{\sigma,b}[\|\hat{\mu} - \mu\|_{L^2}^2]^{1/2} \lesssim N^{-s/(2s+1)}.$$

PROOF. The explicit formula (3.1) for $\mu$ shows that $\|\mu\|_{H^s}$ is uniformly bounded over $\Theta_s$. This implies that the bias term satisfies

$$\|\mu - \pi_J \mu\|_{L^2} \lesssim 2^{-Js}\|\mu\|_{H^s} \sim N^{-s/(2s+1)},$$

uniformly over $\Theta_s$. Since $\hat{\mu}_\lambda$ is an unbiased estimator of $\langle \mu, \psi_\lambda \rangle$, we can apply the variance estimates of Lemma 6.2 to obtain

$$\mathbb{E}_{\sigma,b}[\|\hat{\mu} - \pi_J \mu\|_{L^2}^2] = \sum_{|\lambda| \leq J} \mathrm{Var}_{\sigma,b}[\hat{\mu}_\lambda] \lesssim 2^J N^{-1},$$

which—in combination with the uniformity of the constants involved—gives the announced upper bound. □

4.2. *Spectral approximation.* We shall rely on the spectral approximation results given in Chatelin (1983); compare also Kato (1995). Since for $\hat{\sigma}(\cdot)$ we have to estimate not only the eigenvalue $u_1$, but also its derivative $u_1'$, we will be working in the $L^2$-Sobolev space $H^1$. The general idea is that the error in the eigenvalue and in the eigenfunction can be controlled by the error of the operator on the eigenspace, once the overall error measured in the operator norm is small. Let $R(T, z) = (T - z\,\mathrm{Id})^{-1}$ denote the resolvent map of the operator $T$, $\sigma(T)$ its spectrum and $B(x, r)$ the closed ball of radius $r$ around $x$.

PROPOSITION 4.2. *Suppose a bounded linear operator $T$ on a Hilbert space has a simple eigenvalue $\kappa$ such that $\sigma(T) \cap B(\kappa, \rho) = \{\kappa\}$ holds for some $\rho > 0$. Let $T_\varepsilon$ be a second linear operator with $\|T_\varepsilon - T\| < \frac{R}{2}$, where $R := (\sup_{z \in B(\kappa, \rho)} \|R(T, z)\|)^{-1}$. Then the operator $T_\varepsilon$ has a simple eigenvalue $\kappa_\varepsilon$ in $B(\kappa, \rho)$ and there are eigenvectors $u$ and $u_\varepsilon$ with $Tu = \kappa u$, $T_\varepsilon u_\varepsilon = \kappa_\varepsilon u_\varepsilon$ satisfying*

$$\|u_\varepsilon - u\| \leq \sqrt{8} R^{-1} \|(T_\varepsilon - T)u\|. \tag{4.1}$$

PROOF. We use the resolvent identity and the Cauchy integral representation of the spectral projection $P_\varepsilon$ on the eigenspace of $T_\varepsilon$ contained in $B(\kappa, \rho)$ [see Chatelin (1983), Lemma 6.4]. By the usual Neumann series argument we find formally for an eigenvector $u$ corresponding to $\kappa$,

$$\|u - P_\varepsilon u\| = \frac{1}{2\pi} \left\| \oint_{B(\kappa, \rho)} \frac{R(T_\varepsilon, z)}{\kappa - z} dz \, (T_\varepsilon - T)u \right\|$$

$$\leq \frac{1}{2\pi} 2\pi\rho \sup_{z \in B(\kappa, \rho)} \|R(T_\varepsilon, z)\| \rho^{-1} \|(T_\varepsilon - T)u\|$$

$$\leq \sup_{z \in B(\kappa, \rho)} \frac{\|R(T, z)\|}{1 - \|R(T, z)\| \|T_\varepsilon - T\|} \|(T_\varepsilon - T)u\|$$

$$= (R - \|T_\varepsilon - T\|)^{-1} \|(T_\varepsilon - T)u\|.$$

Hence, for $\|T_\varepsilon - T\| < \frac{R}{2}$ this calculation is a posteriori justified and yields $\|u - P_\varepsilon u\| < 2R^{-1}\|(T_\varepsilon - T)u\|$. Applying $\|(T_\varepsilon - T)u\| < \frac{R}{2}\|u\|$ once again, we see that the projection $P_\varepsilon$ cannot be zero. Consequently there must be a part of the spectrum of $T_\varepsilon$ in $B(\kappa, \rho)$. By the argument in the proof of Theorem 5.22 in Chatelin (1983) this part consists of a simple eigenvalue $\kappa_\varepsilon$.

It remains to find eigenvectors that are close, too. Observe that, for arbitrary Hilbert space elements $g, h$ with $\|g\| = \|h\| = 1$ and $\langle g, h \rangle \geq 0$,
$$\|g - h\|^2 = 2 - 2\langle g, h \rangle \leq 2(1 + \langle g, h \rangle)(1 - \langle g, h \rangle) = 2\|g - \langle g, h \rangle h\|^2$$
holds. We substitute for $g$ and $h$ the normalized eigenvectors $u$ and $u_\varepsilon$ with $\langle u, u_\varepsilon \rangle \geq 0$; note that oblique projections only enlarge the right-hand side and thus infer (4.1). □

COROLLARY 4.3. *Under the conditions of Proposition 4.2 there is a constant $C = C(R, \|T\|)$ such that $|\kappa_\varepsilon - \kappa| \leq C\|(T_\varepsilon - T)u\|$.*

PROOF. The inverse triangle inequality yields
$$|\kappa_\varepsilon - \kappa| = |\|T_\varepsilon u_\varepsilon\| - \|Tu\|| \leq \|T_\varepsilon(u_\varepsilon - u) + (T_\varepsilon - T)u\|$$
$$\leq (\|T\| + \|T_\varepsilon - T\|)\|u_\varepsilon - u\| + \|(T_\varepsilon - T)u\|$$
$$\leq \left(\left(\|T\| + \frac{R}{2}\right)\sqrt{8}R^{-1} + 1\right)\|(T_\varepsilon - T)u\|,$$
where the last line follows from Proposition 4.2. □

4.3. *Bias estimates.* In a first estimation step, we bound the deterministic error due to the finite-dimensional projection $\pi_J^\mu P_\Delta$ of $P_\Delta$. We start with a lemma stating that $\pi_J^\mu$ and $\pi_J$ have similar approximation properties.

LEMMA 4.4. *Let $m : [0, 1] \to [m_0, m_1]$ be a measurable function with $m_1 \geq m_0 > 0$. Denote by $\pi_J^m$ the $L^2(m)$-orthogonal projection onto the multiresolution space $V_J$. Then there is a constant $C = C(m_0, m_1)$ such that*
$$\|(\mathrm{Id} - \pi_J^m)f\|_{H^1} \leq C\|(\mathrm{Id} - \pi_J)f\|_{H^1} \qquad \forall f \in H^1([0, 1]).$$

PROOF. The norm equivalence $m_0\|g\|_{L^2} \leq \|g\|_m \leq m_1\|g\|_{L^2}$ implies
$$\|\pi_J^m\|_{L^2 \to L^2} \leq m_1 m_0^{-1}\|\pi_J^m\|_{L^2(m) \to L^2(m)} = m_1 m_0^{-1}.$$
On the other hand, the Bernstein inequality in $V_J$ and the Jackson inequality for $\mathrm{Id} - \pi_J$ in $H^1$ and $L^2$ yield, for $f \in H^1$,
$$\|(\mathrm{Id} - \pi_J^\mu)f\|_{H^1} = \|(\mathrm{Id} - \pi_J^\mu)(\mathrm{Id} - \pi_J)f\|_{H^1}$$
$$\leq \|(\mathrm{Id} - \pi_J)f\|_{H^1} + \|\pi_J^\mu(\mathrm{Id} - \pi_J)f\|_{H^1}$$
$$\lesssim \|(\mathrm{Id} - \pi_J)f\|_{H^1} + 2^J\|\pi_J^\mu(\mathrm{Id} - \pi_J)f\|_{L^2}$$
$$\lesssim \|(\mathrm{Id} - \pi_J)f\|_{H^1} + \|\pi_J^\mu\|_{L^2 \to L^2}\|(\mathrm{Id} - \pi_J)f\|_{H^1},$$

where the constants depend only on the approximation spaces. □

PROPOSITION 4.5. *Uniformly over $\Theta_s$ we have $\|\pi_J^\mu P_\Delta - P_\Delta\|_{H^1 \to H^1} \lesssim 2^{-Js}$.*

PROOF. The transition density $p_\Delta$ is the kernel of the operator $P_\Delta$. Hence, from Lemma 6.7 it follows that $P_\Delta : H^1 \to H^{s+1}$ is continuous with a uniform norm bound over $\Theta_s$. Lemma 4.4 yields

$$\|(P_\Delta - \pi_J^\mu P_\Delta)f\|_{H^1} \lesssim \|\operatorname{Id} - \pi_J\|_{H^{s+1} \to H^1} \|f\|_{H^1}.$$

The Jackson inequality in $H^1$ gives the result. □

COROLLARY 4.6. *Let $\kappa_1^J$ be the largest eigenvalue smaller than 1 of $\pi_J^\mu$ with eigenfunction $u_1^J$. Then uniformly over $\Theta_s$ the following estimate holds:*

$$|\kappa_1^J - \kappa_1| + \|u_1^J - u_1\|_{H^1} \lesssim 2^{-Js}.$$

PROOF. We are going to apply Proposition 4.2 on the space $H^1$ and its Corollary 4.3. In view of Proposition 4.5, it remains to establish the existence of uniformly strictly positive values for $\rho$ and $R$ over $\Theta_s$. The uniform separation of $\kappa_1$ from the rest of the spectrum is the content of Proposition 6.5.

For the choice of $\rho$ in Proposition 6.5 we search a uniform bound $R$. If we regard $P_\Delta$ on $L^2(\mu)$, then $P_\Delta$ is self-adjoint and satisfies $\|R(P_\Delta, z)\| = \operatorname{dist}(z, \sigma(P_\Delta))^{-1}$ [see Chatelin (1983), Proposition 2.32].

By Lemma 6.3 and the commutativity between $P_\Delta$ and $L$ we conclude

$$\|R(P_\Delta, z)f\|_{H^1} \sim \|(\operatorname{Id} - L)^{1/2} R(P_\Delta, z)f\|_\mu$$
$$\leq \|R(P_\Delta, z)\| \|(\operatorname{Id} - L)^{1/2} f\|_\mu$$
$$\sim \operatorname{dist}(z, \sigma(P_\Delta))^{-1} \|f\|_{H^1}.$$

Hence, $\|R(P_\Delta, z)\|_{H^1 \to H^1} \lesssim \rho^{-1}$ holds uniformly over $z \in B(\kappa, \rho)$ and $(\sigma, b) \in \Theta_s$. □

REMARK 4.7. The Kato–Temple inequality [Chatelin (1983), Theorem 6.21] on $L^2(\mu)$ even establishes the so-called superconvergence $|\kappa_1^J - \kappa_1| \lesssim 2^{-2Js}$.

4.4. *Variance estimates.* To bound the stochastic error on the finite-dimensional space $V_J$, we return to vector-matrix notation and look for a bound on the error involved in the estimators $\hat{\mathbf{G}}$ and $\hat{\mathbf{P}}_\Delta$. The Euclidean norm is denoted by $\|\cdot\|_{l^2}$.

LEMMA 4.8. *For any vector* $\mathbf{v} \in \mathbb{R}^{|V_J|}$ *we have, uniformly over* $\Theta_s$,
$$\mathbb{E}_{\sigma,b}[\|(\hat{\mathbf{G}} - \mathbf{G})\mathbf{v}\|_{l^2}^2] \lesssim \|\mathbf{v}\|_{l^2}^2 N^{-1} 2^J.$$

PROOF. We obtain, by (6.1) in Lemma 6.2,

$$\mathbb{E}_{\sigma,b}[\|(\mathbf{G} - \hat{\mathbf{G}})\mathbf{v}\|_{l^2}^2]$$
$$= \sum_{|\lambda| \leq J} \mathbb{E}_{\sigma,b}\left[\frac{1}{N^2}\left(\mathbb{E}_{\sigma,b}[\psi_\lambda(X_0)v(X_0)] - \frac{1}{2}\psi_\lambda(X_0)v(X_0)\right.\right.$$
$$\left.\left. - \frac{1}{2}\psi_\lambda(X_{N\Delta})v(X_{N\Delta}) - \sum_{n=1}^{N-1} \psi_\lambda(X_{n\Delta})v(X_{n\Delta})\right)^2\right]$$
$$\lesssim \sum_{|\lambda| \leq J} N^{-1} \mathbb{E}_{\sigma,b}[(\psi_\lambda(X_0)v(X_0))^2]$$
$$\lesssim N^{-1}\left\|\sum_{|\lambda| \leq J} \psi_\lambda^2\right\|_\infty \|v\|_{L^2}^2 \|\mu\|_\infty \lesssim N^{-1} 2^J \|\mathbf{v}\|_{l^2}^2,$$

as asserted. □

LEMMA 4.9. *For any vector* $\mathbf{v}$ *we have, uniformly over* $\Theta_s$,
$$\mathbb{E}_{\sigma,b}[\|(\hat{\mathbf{P}}_\Delta - \mathbf{P}_\Delta)\mathbf{v}\|_{l^2}^2] \lesssim \|\mathbf{v}\|_{l^2}^2 N^{-1} 2^J.$$

PROOF. We obtain, by (6.2) in Lemma 6.2,

$$\mathbb{E}_{\sigma,b}[\|(\hat{\mathbf{P}}_\Delta - \mathbf{P}_\Delta)v\|_{l^2}^2]$$
$$= \sum_{|\lambda| \leq J} \mathbb{E}_{\sigma,b}\left[\left(\frac{1}{N}\sum_{n=1}^N \psi_\lambda(X_{(n-1)\Delta})v(X_{n\Delta}) - \mathbb{E}_{\sigma,b}[\psi_\lambda(X_0)v(X_\Delta)]\right)^2\right]$$
$$\lesssim \sum_{|\lambda| \leq J} N^{-1} \mathbb{E}_{\sigma,b}[(\psi_\lambda(X_0)v(X_1))^2]$$
$$\leq N^{-1} \sum_{|\lambda| \leq J} \|\psi_\lambda\|_{L^2}^2 \|v\|_{L^2}^2 \|\mu p_\Delta\|_\infty \lesssim N^{-1} 2^J \|\mathbf{v}\|_{l^2}^2,$$

as asserted. □

DEFINITION 4.10. We introduce the random set
$$\mathcal{R} = \mathcal{R}_{J,N} := \{\|\hat{\mathbf{G}} - \mathbf{G}\| \leq \tfrac{1}{2}\|\mathbf{G}^{-1}\|^{-1}\}.$$

REMARK 4.11. Since $\mathbf{G}$ is invertible, so is $\hat{\mathbf{G}}$ on $\mathcal{R}$ with $\|\hat{\mathbf{G}}^{-1}\| \leq 2\|\mathbf{G}^{-1}\|$ by the usual Neumann series argument.

LEMMA 4.12. *Uniformly over $\Theta_s$ we have $\mathbb{P}_{\sigma,b}(\Omega \setminus \mathcal{R}) \lesssim N^{-1}2^{2J}$.*

PROOF. By the classical Hilbert–Schmidt norm inequality,
$$\|\hat{\mathbf{G}} - \mathbf{G}\|_{l^2 \to l^2}^2 \leq \sum_{|\lambda| \leq J} \|(\hat{\mathbf{G}} - \mathbf{G})\mathbf{e}_\lambda\|_{l^2 \to l^2}^2$$
holds with unit vectors $(\mathbf{e}_\lambda)$. Then Lemma 4.8 gives uniformly
$$\mathbb{E}_{\sigma,b}[\|\hat{\mathbf{G}} - \mathbf{G}\|_{l^2 \to l^2}^2] \lesssim N^{-1}2^{2J}.$$
Since the spaces $L^2([0,1])$ and $L^2(\mu)$ are isomorphic with uniform isomorphism constants, $\|\mathbf{G}^{-1}\| \sim 1$ holds uniformly over $\Theta_s$ and the assertion follows from Chebyshev's inequality. □

PROPOSITION 4.13. *For any $\varepsilon > 0$ we have, uniformly over $\Theta_s$,*
$$\mathbb{P}_{\sigma,b}(\mathcal{R} \cap \{\|\hat{\mathbf{G}}^{-1}\hat{\mathbf{P}}_\Delta - \mathbf{G}^{-1}\mathbf{P}_\Delta\| \geq \varepsilon\}) \lesssim N^{-1}2^{2J}\varepsilon^{-2}.$$

PROOF. First, we separate the different error terms:
$$\hat{\mathbf{G}}^{-1}\hat{\mathbf{P}}_\Delta - \mathbf{G}^{-1}\mathbf{P}_\Delta = \hat{\mathbf{G}}^{-1}(\hat{\mathbf{P}}_\Delta - \mathbf{P}_\Delta) + (\hat{\mathbf{G}}^{-1} - \mathbf{G}^{-1})\mathbf{P}_\Delta$$
$$= \hat{\mathbf{G}}^{-1}((\hat{\mathbf{P}}_\Delta - \mathbf{P}_\Delta) + (\mathbf{G} - \hat{\mathbf{G}})\mathbf{G}^{-1}\mathbf{P}_\Delta).$$
On the set $\mathcal{R}$ we obtain, by Remark 4.11,
$$\|\hat{\mathbf{G}}^{-1}\hat{\mathbf{P}}_\Delta - \mathbf{G}^{-1}\mathbf{P}_\Delta\| \leq \|\hat{\mathbf{G}^{-1}}\|(\|\hat{\mathbf{P}}_\Delta - \mathbf{P}_\Delta\| + \|\mathbf{G} - \hat{\mathbf{G}}\|\|\mathbf{G}^{-1}\|\|\mathbf{P}_\Delta\|)$$
$$\leq 2\|\mathbf{G}^{-1}\|(\|\hat{\mathbf{P}}_\Delta - \mathbf{P}_\Delta\| + \|\mathbf{G} - \hat{\mathbf{G}}\|\|\mathbf{G}^{-1}\|\|\mathbf{P}_\Delta\|)$$
$$\lesssim \|\hat{\mathbf{P}}_\Delta - \mathbf{P}_\Delta\| + \|\mathbf{G} - \hat{\mathbf{G}}\|.$$
By Lemmas 4.8 and 4.9 and the Hilbert–Schmidt norm estimate (cf. the proof of Lemma 4.12) we obtain the uniform norm bound over $\Theta_s$,
$$\mathbb{E}_{\sigma,b}[\|\hat{\mathbf{G}}^{-1}\hat{\mathbf{P}}_\Delta - \mathbf{G}^{-1}\mathbf{P}_\Delta\|^2 \mathbf{1}_\mathcal{R}] \lesssim N^{-1}2^{2J}.$$
It remains to apply Chebyshev's inequality. □

Having established the weak consistency of the estimators in matrix norm, we now bound the error on the eigenspace.

PROPOSITION 4.14. *Let $\mathbf{u}_1^J$ be the vector associated with the normalized eigenfunction $u_1^J$ of $\pi_J^\mu P_\Delta$ with eigenvalue $\kappa_1^J$. Then uniformly over $\Theta_s$ the following risk bound holds:*
$$\mathbb{E}_{\sigma,b}[\|(\hat{\mathbf{G}}^{-1}\hat{\mathbf{P}}_\Delta - \mathbf{G}^{-1}\mathbf{P}_\Delta^J)\mathbf{u}_1^J\|_{l^2}^2 \mathbf{1}_\mathcal{R}] \lesssim N^{-1}2^J.$$

PROOF. By the same separation of the error terms on $\mathcal{R}$ as in the preceding proof and by Lemmas 4.8 and 4.9 we find

$$\mathbb{E}_{\sigma,b}[\|(\hat{\mathbf{G}}^{-1}\hat{\mathbf{P}}_\Delta - \mathbf{G}^{-1}\mathbf{P}_\Delta^J)\mathbf{u}_1^J\|_{l^2}^2 \mathbf{1}_{\mathcal{R}}]$$
$$\leq 8\|\mathbf{G}^{-1}\|^2(\mathbb{E}[\|(\hat{\mathbf{P}}_\Delta - \mathbf{P}_\Delta^J)\mathbf{u}_1^J\|_{l^2}^2] + \mathbb{E}[\|(\mathbf{G}-\hat{\mathbf{G}})\kappa_1^J \mathbf{u}_1^J\|_{l^2}^2]) \lesssim N^{-1}2^J.$$

The uniformity over $\Theta_s$ follows from the respective statements in the lemmas. $\square$

COROLLARY 4.15. *Let $\hat{\kappa}_1$ be the second largest eigenvalue of the matrix $\hat{\mathbf{G}}^{-1}\hat{\mathbf{P}}_\Delta$ with eigenvector $\hat{\mathbf{u}}_1$. If $\hat{\mathbf{G}}$ is not invertible or if $\|\hat{\mathbf{u}}_1\|_{l^2} \geq 2\sup_{\Theta_s}\|u_1\|_{L^2}$ holds, put $\hat{\mathbf{u}}_1 := \mathbf{0}$, $\hat{\kappa}_1 := 0$. If $N^{-1}2^{2J} \to 0$ holds, then uniformly over $\Theta_s$ the following bounds hold for $N, J \to \infty$:*

(4.2) $$\mathbb{E}_{\sigma,b}[(|\hat{\kappa}_1 - \kappa_1^J|^2 + \|\hat{\mathbf{u}}_1 - \mathbf{u}_1^J\|_{l^2}^2)\mathbf{1}_{\mathcal{R}}] \lesssim N^{-1}2^J,$$

(4.3) $$\mathbb{E}_{\sigma,b}[\|\hat{u}_1 - u_1^J\|_{H^1}^2] \lesssim N^{-1}2^{3J}.$$

PROOF. For the proof of (4.2) we apply Proposition 4.2 using the Euclidean Hilbert space $\mathbb{R}^{V_J}$ and Corollary 4.3. Then Proposition 4.14 in connection with Proposition 4.13 (using $\varepsilon < R/2$ and $N^{-1}2^{2J} \to 0$) yields the correct asymptotic rate on the event $\mathcal{R}$. For the uniform choice of $\rho$ and $R$ for $\mathbf{G}^{-1}\mathbf{P}_\Delta^J$ in Proposition 4.2 just use the corresponding result for $P_\Delta$ and the convergence $\|\pi_J^\mu P_\Delta - P_\Delta\| \to 0$.

The precaution taken for undefined or too large $\hat{\mathbf{u}}_1$ is necessary for the event $\Omega \setminus \mathcal{R}$. Since the estimators $\hat{\kappa}_1$ and $\hat{\mathbf{u}}_1$ are now kept artificially bounded, the rate $\mathbb{P}_{\sigma,b}(\Omega \setminus \mathcal{R}) \lesssim N^{-1}2^{2J}$ established in Lemma 4.12 suffices to bound the risk on $\Omega \setminus \mathcal{R}$. Hence, the second estimate (4.3) is a consequence of (4.2) and the Bernstein inequality $\|\hat{u}_1 - u_1^J\|_{H^1} \lesssim 2^J \|\hat{u}_1 - u_1^J\|_{L^2}$. $\square$

REMARK 4.16. The main result of this section, namely (4.3), can be extended to $p$th moments for all $p \in (1, \infty)$:

(4.4) $$\mathbb{E}_{\sigma,b}[\|\hat{u}_1 - u_1^J\|_{H^1}^p]^{1/p} \lesssim N^{-1/2}2^{3J/2}.$$

Indeed, tracing back the steps, it suffices to obtain bounds on the moments of order $p$ in Lemmas 4.8 and 4.9, which on their part rely on the mixing statement in Lemma 6.2. By arguments based on the Riesz convexity theorem this last lemma generalizes to the corresponding bounds for $p$th moments, as derived in Section VII.4 of Rosenblatt (1971). For the sake of clarity we have restricted ourselves to the case $p = 2$ here.

4.5. *Upper bound for $\sigma(\cdot)$.* By Corollary 4.15 and our choice of $J$, $2^J \sim N^{1/(2s+3)}$,

$$\sup_{(\sigma,b) \in \Theta_s} \mathbb{E}_{\sigma,b}[|\hat{\kappa}_1 - \kappa_1^J|^2 + \|\hat{u}_1 - u_1^J\|_{H^1}^2] \lesssim N^{-1} 2^{3J} \sim N^{-2s/(2s+3)}$$

holds. Using this estimate and the estimate for $\hat{\mu}$ in Proposition 4.1, the risk of the plug-in estimator $\hat{\sigma}^2(\cdot)$ in (3.12) is bounded as asserted in the theorem. We only have to ensure that the stochastic error does not increase from the plug-in and that the denominator is uniformly bounded away from zero. Using the Cauchy–Schwarz inequality and Remark 4.16 on the higher moments of our estimators, we encounter no problem in the first case. The second issue is dealt with by using the lower bound $c_{a,b} > 0$ in Proposition 6.5 so that an improvement of the estimate for the denominator by using

$$\widehat{\mu u}_1 := \max(\hat{\mu}\hat{u}_1, c_{a,b}) \tag{4.5}$$

instead of $\hat{\mu}\hat{u}_1$ guarantees the uniform lower bound away from zero.

4.6. *Upper bound for $b(\cdot)$.* Since $b(\cdot) = S'(\cdot)/\mu(\cdot)$ holds, it suffices to discuss how to estimate $S'(\cdot)$, which amounts to estimating the eigenfunction $u_1$ in $H^2$-norm; compare with (3.6). Substituting $H^2$ for $H^1$ in Proposition 4.5 and its proof, we obtain the bound

$$\|\pi_J^\mu P_\Delta - P_\Delta\|_{H^2 \to H^2} \lesssim 2^{-J(s-1)},$$

because $\|\mathrm{Id} - \pi_J\|_{H^{s+1} \to H^2}$ is of this order. As in Corollary 4.6 this is also the rate for the bias estimate. The only fine point is the uniform norm equivalence $\|f\|_{H^2} \sim \|(\mathrm{Id} - L)f\|_\mu$ for $f \in \mathrm{dom}(L)$, which follows by the methodology of perturbation and similarity arguments given in Section VI.4b of Engel and Nagel (2000). We omit further details.

The variance estimate is exactly the same. From (4.2) we infer, by Bernstein's inequality for $H^2$ and the estimate of $\mathbb{P}_{\sigma,b}(\Omega \setminus \mathcal{R})$,

$$\mathbb{E}_{\sigma,b}[\|\hat{u}_1 - u_1^J\|_{H^2}^2] \lesssim N^{-1} 2^{5J}.$$

Therefore balancing the bias and variance part of the risk by the choice $2^J \sim N^{1/(2s+3)}$—as before—yields the asserted rate of convergence $N^{-(s-1)/(2s+3)}$.

**5. Proof of the lower bounds.** First, the usual Bayes prior technique is applied for the reduction to a problem of bounding certain likelihood ratios. Then the problem is reduced to that of bounding the $L^2$-distance between the transition probabilities, which is finally accomplished using Hilbert–Schmidt norm estimates and the explicit form of the inverse of the generator.

5.1. *The least favorable prior.* The idea is to perturb the drift and diffusion coefficients of the reflected Brownian motion in such a way that the invariant measure remains unchanged. Let us assume that $\psi$ is a compactly supported wavelet in $H^s$ with one vanishing moment. Without loss of generality we suppose $C > 1 > c > 0$ such that $(1,0) \in \Theta_s$ holds. We put $\psi_{jk} = 2^{j/2}\psi(2^j \cdot - k)$ and denote by $K_j \subset \mathbb{Z}$ a maximal set of indices $k$ such that $\mathrm{supp}(\psi_{jk}) \subset [a,b]$ and $\mathrm{supp}(\psi_{jk}) \cap \mathrm{supp}(\psi_{jk'}) = \varnothing$ holds for all $k, k' \in K_j$, $k \neq k'$. Furthermore, we set $\gamma \sim 2^{-j(s+1/2)}$ such that, for all $\varepsilon = \varepsilon(j) \in \{-1, +1\}^{|K_j|}$,

$$(\sqrt{2S_\varepsilon}, (S_\varepsilon)') \in \Theta_s \qquad \text{with } S_\varepsilon(x) := S_\varepsilon(j, x) := \left(2 + \gamma \sum_{k \in K_j} \varepsilon_k \psi_{jk}(x)\right)^{-1}.$$

We consider the corresponding diffusions with generator

$$L_{S_\varepsilon} f(x) := (S_\varepsilon f')'(x) := S_\varepsilon(x) f''(x) + S'_\varepsilon(x) f'(x), \qquad f \in \mathrm{dom}(L).$$

Hence, the invariant measure is the Lebesgue measure on $[0, 1]$.

The usual Assouad cube techniques [e.g., see Korostelev and Tsybakov (1993)] give, for any estimator $\hat{\sigma}(\cdot)$ and for $N \in \mathbb{N}$, $\rho > 0$, the lower bounds

$$(5.1) \qquad \sup_{(\sigma, b) \in \Theta_s} \mathbb{E}_{\sigma, b}[\|\hat{\sigma}^2 - \sigma^2\|^2_{L^2([a,b])}] \gtrsim 2^j \delta^2_{\sigma^2} e^{-\rho} p_0,$$

$$(5.2) \qquad \sup_{(\sigma, b) \in \Theta_s} \mathbb{E}_{\sigma, b}[\|\hat{b} - b\|^2_{L^2([a,b])}] \gtrsim 2^j \delta^2_b e^{-\rho} p_0,$$

where, for all $\varepsilon, \varepsilon'$ with $\|\varepsilon - \varepsilon'\|_{l^1} = 2$ and with $\mathbb{P}_\varepsilon := \mathbb{P}_{\sqrt{2S_\varepsilon}, S'_\varepsilon}$,

$$\delta_{\sigma^2} \leq \|2S_{\varepsilon'} - 2S_\varepsilon\|_{L^2}, \qquad \delta_b \leq \|S'_{\varepsilon'} - S'_\varepsilon\|_{L^2}, \qquad p_0 \leq \mathbb{P}_\varepsilon\left(\left.\frac{d\mathbb{P}_{\varepsilon'}}{d\mathbb{P}_\varepsilon}\right|_{\mathcal{F}_N} > e^{-\rho}\right).$$

We choose $\delta_{\sigma^2} \sim \gamma$ since, for $x \in \mathrm{supp}(\psi_{jk})$ with $\varepsilon_k = -\varepsilon'_k$,

$$S_{\varepsilon'}(x) - S_\varepsilon(x) = \pm 2\gamma \psi_{jk}(x) S_{\varepsilon'}(x) S_\varepsilon(x)$$

and $S_\varepsilon, S_{\varepsilon'} \to \frac{1}{2}$ holds uniformly so that the $L^2$-norm is indeed of order $\gamma$. Equivalently, we find $\delta_b \sim 2^j \gamma$. Due to $\gamma \sim 2^{-j(s+1/2)}$, the proof of the theorem is accomplished, once we have shown that a strictly positive $p_0$ can be chosen for fixed $\rho > 0$ and the asymptotics $2^j \sim N^{1/(2s+3)}$.

5.2. *Reduction to the convergence of transition densities.* If we denote the transition probability densities $\mathbb{P}_\varepsilon(X_\Delta \in dy | X_0 = x)$ by $p_\varepsilon(x, y) \, dy$ and the transition density of reflected Brownian motion by $p_{\mathrm{BM}}$, then we infer, from Proposition 6.4,

$$\lim_{j \to \infty} \sup_{\varepsilon \in \{-1,+1\}^{K_j}} \|p_\varepsilon - p_{\mathrm{BM}}\|_\infty = 0$$

due to $\|S_\varepsilon - \tfrac{1}{2}\|_{C^1} \sim \gamma 2^{3j/2} \to 0$ for $s > 1$. We are now going to use the estimate $-\log(1+x) \leq x^2 - x$, which is valid for all $x \geq -\tfrac{1}{2}$. For $j$ so large that $\|1 - \tfrac{p_{\varepsilon'}}{p_\varepsilon}\|_\infty \leq \tfrac{1}{2}$ holds, the Kullback–Leibler distance can be bounded from above (note that the invariant measure is Lebesgue measure):

$$\mathbb{E}_\varepsilon\left[-\log\left(\frac{d\mathbb{P}_{\varepsilon'}}{d\mathbb{P}_\varepsilon}\bigg|_{\mathcal{F}_N}\right)\right]$$

$$= -\sum_{k=1}^N \mathbb{E}_\varepsilon\left[\log\left(\frac{p_{\varepsilon'}(X_{k-1}, X_k)}{p_\varepsilon(X_{k-1}, X_k)}\right)\right]$$

$$= -N \int_0^1 \int_0^1 \log\left(\frac{p_{\varepsilon'}(x,y)}{p_\varepsilon(x,y)}\right) p_\varepsilon(x,y)\, dy\, dx$$

$$\leq N \int_0^1 \int_0^1 \frac{(p_{\varepsilon'}(x,y) - p_\varepsilon(x,y))^2}{p_\varepsilon(x,y)} - (p_{\varepsilon'}(x,y) - p_\varepsilon(x,y))\, dy\, dx$$

$$= N \int_0^1 \int_0^1 \frac{(p_{\varepsilon'}(x,y) - p_\varepsilon(x,y))^2}{p_\varepsilon(x,y)}\, dy\, dx$$

$$\leq N \|p_\varepsilon^{-1}\|_\infty \|p_{\varepsilon'} - p_\varepsilon\|_{L^2([0,1]^2)}^2.$$

The square root of the Kullback–Leibler distance bounds the total variation distance in order, which by the Chebyshev inequality yields

$$\mathbb{P}_\varepsilon\left(\frac{d\mathbb{P}_{\varepsilon'}}{d\mathbb{P}_\varepsilon}\bigg|_{\mathcal{F}_N} > e^{-\rho}\right) = 1 - \mathbb{P}_\varepsilon\left(\frac{d\mathbb{P}_{\varepsilon'}}{d\mathbb{P}_\varepsilon}\bigg|_{\mathcal{F}_N} - 1 \leq e^{-\rho} - 1\right)$$

$$\geq 1 - \mathbb{E}_\varepsilon\left[\left|\frac{d\mathbb{P}_{\varepsilon'}}{d\mathbb{P}_\varepsilon}\bigg|_{\mathcal{F}_N} - 1\right|\right](1 - e^{-\rho})^{-1}$$

$$= 1 - (1 - e^{-\rho})^{-1} \|(\mathbb{P}_{\varepsilon'} - \mathbb{P}_\varepsilon)|_{\mathcal{F}_N}\|_{\mathrm{TV}}$$

$$\geq 1 - C N^{1/2} \|p_{\varepsilon'} - p_\varepsilon\|_{L^2([0,1]^2)},$$

where $C > 0$ is some constant independent of $\gamma$, $N$, $\varepsilon$ and $j$. Summarizing, we need the estimate

(5.3) $\quad\limsup_{N,j\to\infty} N^{1/2} \|p_{\varepsilon'} - p_\varepsilon\|_{L^2([0,1]^2)} < C^{-1} \quad$ for $2^j \sim N^{1/(2s+3)}$.

5.3. *Convergence of the transition densities.* Observe first that $\|p_{\varepsilon'} - p_\varepsilon\|_{L^2([0,1]^2)}$ is exactly the Hilbert–Schmidt norm distance $\|P_\Delta^{\varepsilon'} - P_\Delta^\varepsilon\|_{\mathrm{HS}}$ between the transition operators derived from $L_{S_{\varepsilon'}}$ and $L_{S_\varepsilon}$ acting on the Hilbert space $L^2([0,1])$. If we introduce

$$V := \left\{f \in L^2([0,1]) \,\bigg|\, \int f = 0\right\} \quad \text{and} \quad V^\perp := \{f \in L^2([0,1]) | f \text{ constant}\},$$

then the transition operators coincide on $V^\perp$ and leave the space $V$ invariant so that $\|P_\Delta^{\varepsilon'} - P_\Delta^\varepsilon\|_{\text{HS}} = \|(P_\Delta^{\varepsilon'} - P_\Delta^\varepsilon)|_V\|_{\text{HS}}$.

We take advantage of the key result that for Lipschitz functions $f$ with Lipschitz constant $\Lambda$ on the union of the spectra of two self-adjoint bounded operators $T_1$ and $T_2$ the continuous functional calculus satisfies

$$(5.4) \qquad \|f(T_1) - f(T_2)\|_{\text{HS}} \leq \Lambda \|T_1 - T_2\|_{\text{HS}};$$

see Kittaneh (1985). We proceed by bounding the Hilbert–Schmidt norm of the difference of the inverses of the generators and by then transferring this bound to the transition operators via (5.4). By the functional calculus for operators on $V$, the function $f(z) = \exp(\Delta(z^{-1}))$ sends $(L_\varepsilon|_V)^{-1}$ to $P_\Delta^\varepsilon|_V$. Moreover, $f$ is uniformly Lipschitz continuous on $(-\infty, 0)$ due to $\Lambda := \sup_{z<0} |f'(z)| = 4\Delta^{-1} e^{-2} < \infty$. Thus, we arrive at

$$\|p_{\varepsilon'} - p_\varepsilon\|_{L^2([0,1]^2)} = \|(P_\Delta^{\varepsilon'} - P_\Delta^\varepsilon)|_V\|_{\text{HS}} \leq \Lambda \|(L_{\varepsilon'}|_V)^{-1} - (L_\varepsilon|_V)^{-1}\|_{\text{HS}}.$$

The inverse of the generator $L_\varepsilon$ on $V$ has for $g \in V$ the explicit form

$$(5.5) \qquad (L_\varepsilon|_V)^{-1} g(x) = \int_0^1 \left( \int_y^1 S_\varepsilon^{-1}(v)(v - \mathbf{1}_{[x,1]}(v)) \, dv \right) g(y) \, dy.$$

Using $|S_{\varepsilon'}^{-1} - S_\varepsilon^{-1}| = 2\gamma \psi_{jk}$ for some $k \in K_j$ and denoting by $\Psi$ the primitive of $\psi$ with compact support, we obtain

$$\|(L_{\varepsilon'}|_V)^{-1} - (L_\varepsilon|_V)^{-1}\|_{\text{HS}}^2$$

$$= \int_0^1 \int_0^1 \left( \int_y^1 2\gamma \psi_{jk}(v)(v - \mathbf{1}_{[x,1]}(v)) \, dv \right)^2 dx \, dy$$

$$= 4\gamma^2 2^{-j} \int_0^1 \int_0^1 \left( -\Psi(2^j y - k)y \right.$$

$$\left. - \int_y^1 \Psi(2^j v - k) \, dv + \Psi(2^j(x \vee y) - k) \right)^2 dx \, dy$$

$$\lesssim \gamma^2 2^{-j} \|\Psi(2^j \cdot)\|_{L^2}^2 \sim \gamma^2 2^{-2j}.$$

Consequently, $\|p_{\varepsilon'} - p_\varepsilon\|_{L^2}^2 \sim 2^{-j(2s+3)}$ holds with an arbitrarily small constant if $\gamma 2^{j(s+1/2)}$ is chosen sufficiently small. Hence, the estimate (5.3) is valid for this choice and the asymptotics $N 2^{-j(2s+3)} \to 1$, which remained to be proved.

**6. Technical results.** We shall need several technical results, mainly to describe the dependence of certain quantities on the underlying diffusion parameters. The following result is in close analogy with Section IV.5 in Bass (1998).

LEMMA 6.1. *The second largest eigenvalue $\nu_1$ of the infinitesimal generator $L_{\sigma,b}$ can be bounded away from zero:*

$$\nu_1 \leq -\inf_{x \in [0,1]} S(x) =: -s_0.$$

*This eigenvalue is simple and the corresponding eigenfunction $f_1$ is monotone.*

PROOF. The variational characterization of $\nu_1$ [Davies (1995), Section 4.5] and partial integration yield

$$\nu_1 = \sup_{\substack{\|f\|_\mu = 1 \\ \langle f, \mathbf{1} \rangle_\mu = 0}} \langle Lf, f \rangle_\mu = -\inf_{\substack{\|f\|_\mu = 1 \\ \langle f, \mathbf{1} \rangle_\mu = 0}} \int_0^1 S(x) f'(x)^2 \, dx.$$

Given the derivative $f'$, the function $f \in \text{dom}(L)$ with $\langle f, \mathbf{1} \rangle_\mu = 0$ is uniquely determined. Setting $M(x) := \mu([0, x])$, this function $f$ satisfies

$$0 = \int_0^1 \left( f(0) + \int_0^x f'(y) \, dy \right) \mu(x) \, dx = f(0) + \int_0^1 f'(y)(1 - M(y)) \, dy.$$

For $f, g \in L^2(\mu)$ with $\langle f, \mathbf{1} \rangle_\mu = \langle g, \mathbf{1} \rangle_\mu = 0$ we find

$$\langle f, g \rangle_\mu = \int_0^1 \left( f(0) + \int_0^x f'(y) \, dy \right) \left( g(0) + \int_0^x g'(z) \, dz \right) \mu(x) \, dx$$

$$= f(0)g(0) - f(0) \int_0^1 g'(z)(1 - M(z)) \, dz$$

$$\quad - g(0) \int_0^1 f'(y)(1 - M(y)) \, dy$$

$$\quad + \int_0^1 \int_0^1 f'(y) g'(z)(1 - M(y \vee z)) \, dy \, dz$$

$$= \int_0^1 \int_0^1 (M(y \wedge z) - M(y) M(z)) f'(y) g'(z) \, dy \, dz$$

$$=: \int_0^1 \int_0^1 m(y, z) f'(y) g'(z) \, dy \, dz.$$

The kernel $m(y, z)$ is positive on $(0, 1)^2$ and bounded by 1, whence we obtain, by regarding $u = f'$,

$$-\nu_1 = \inf_u \frac{\int_0^1 S(x) u(x)^2 \, dx}{\int_0^1 \int_0^1 m(y, z) u(y) u(z) \, dy \, dz} \geq \frac{s_0 \|u\|_{L^2}^2}{\|u\|_{L^1}^2} \geq s_0.$$

If the derivative of an eigenfunction $f_1$ changed sign, we could write $f_1' = u^+ - u^-$ with two nonnegative functions $u^+, u^-$ that are nontrivial. However, this would entail that the antiderivative $f_0$ of $f_0' := u^+ + u^-$ satisfies $\langle Lf_0, f_0\rangle_\mu = \langle Lf_1, f_1\rangle_\mu$, while $\|f_0\|_\mu$ would be strictly greater than $\|f_1\|_\mu$ due to the positivity of $m(u,v)$. This contradicts the variational characterization of $\nu_1$ so that all eigenfunctions corresponding to $\nu_1$ are monotone. Consequently, for any two eigenfunctions $f_1$ and $g_1$ the integrand in

$$\langle f_1, g_1\rangle_\mu = \int_0^1 \int_0^1 m(y,z) f_1'(y) g_1'(z)\, dy\, dz$$

does not change sign and the whole integral does not vanish. We infer that the eigenspace of $\nu_1$ cannot contain two orthogonal functions and is thus one-dimensional. $\square$

LEMMA 6.2. *For $H_1, H_2 \in L^2([0,1])$ we have the following two uniform variance estimates over $\Theta_s$:*

(6.1) $$\operatorname{Var}_{\sigma,b}\left[\frac{1}{N}\sum_{n=1}^N H_1(X_{n\Delta})\right] \lesssim N^{-1}\mathbb{E}_{\sigma,b}[H_1(X_0)^2],$$

(6.2) $$\operatorname{Var}_{\sigma,b}\left[\frac{1}{N}\sum_{n=1}^N H_1(X_{(n-1)\Delta})H_2(X_{n\Delta})\right] \lesssim N^{-1}\mathbb{E}_{\sigma,b}[H_1(X_0)^2 H_2(X_1)^2].$$

PROOF. Due to the uniform spectral gap $s_0$ over $\Theta_s$ (Proposition 6.1), $P_\Delta$ satisfies $\|P_\Delta f\|_\mu \leq \gamma \|f\|_\mu$ with $\gamma := e^{-|s_0|\Delta} < 1$ for all $f \in L^2(\mu)$ with $\langle f, \mathbf{1}\rangle_\mu = 0$.

We obtain the first estimate by considering the centered random variables $f_1(X_{k\Delta}) := H_1(X_{k\Delta}) - \mathbb{E}_{\sigma,b}[H_1(X_{k\Delta})]$, $k \in \mathbb{N}$:

$$\operatorname{Var}_{\sigma,b}\left[\sum_{n=1}^N H_1(X_{n\Delta})\right] = \sum_{m,n=1}^N \mathbb{E}_{\sigma,b}[f_1(X_{m\Delta}) f_1(X_{n\Delta})] = \sum_{m,n=1}^N \langle f_1, P_\Delta^{|m-n|} f_1\rangle_\mu$$

$$\leq \sum_{m,n=1}^N \gamma^{|m-n|} \|f_1\|_\mu^2 \lesssim N \mathbb{E}_{\sigma,b}[H_1(X_0)^2].$$

The second estimate follows along the same lines. Merely observe that for $m > n$, by the projection property of conditional expectations

$$\mathbb{E}_{\sigma,b}[H_1(X_{(n-1)\Delta}) H_2(X_{n\Delta}) H_1(X_{(m-1)\Delta}) H_2(X_{m\Delta})]$$
$$= \langle H_1 \cdot (P_\Delta H_2), P_\Delta^{m-n-1}(H_1 \cdot (P_\Delta H_2))\rangle_\mu$$

holds, where "·" is the usual multiplication operator. $\square$

LEMMA 6.3. *Uniformly over $\Theta_s$ the following norm equivalence holds:*

$$\|f\|_{H^1} \sim \|(\operatorname{Id} - L)^{1/2} f\|_\mu \qquad \text{for all } f \in H^1.$$

PROOF. The invariant measure $\mu$ and the function $S$ are uniformly bounded away from zero and infinity so that we obtain, with uniform constants for $f \in \mathrm{dom}(L)$,

$$\|f\|_{H^1}^2 = \langle f, f \rangle + \langle f', f' \rangle \sim \langle f, f \rangle_\mu + \langle S f', f' \rangle$$
$$= \langle (\mathrm{Id} - L) f, f \rangle_\mu = \|(\mathrm{Id} - L)^{1/2} f\|_\mu^2.$$

By an approximation argument this extends to all $f \in H^1 = \mathrm{dom}(L^{1/2})$. □

PROPOSITION 6.4. *Suppose $((\sigma_n(\cdot), b_n(\cdot)) \in \Theta_s$, $n \geq 0$, and*

$$\lim_{n \to \infty} \|\sigma_n - \sigma_0\|_\infty = 0, \qquad \lim_{n \to \infty} \|b_n - b_0\|_\infty = 0.$$

*Then the corresponding transition probabilities $(p_t^{(n)})$ converge uniformly:*

$$\lim_{n \to \infty} \|p_t^{(n)}(\cdot, \cdot) - p_t^{(0)}(\cdot, \cdot)\|_\infty = 0.$$

PROOF. An application of the results by Stroock and Varadhan (1971) yields that the corresponding diffusion processes $X^{(n)}$ converge weakly to $X^{(0)}$ for any fixed initial value $X^{(n)}(0) = x$. This implies in particular

$$\lim_{n \to \infty} \int_0^1 p_t^{(n)}(x, y) \varphi(y)\, dy = \int_0^1 p_t^{(0)}(x, y) \varphi(y)\, dy$$

for all test functions $\varphi \in L^\infty([0, 1])$ and all $x \in [0, 1]$.

On the other hand, the functions $(p_t^{(n)})_n$ form a relatively compact subset of $\mathcal{C}([0, 1]^2)$ by Proposition 6.7 and Sobolev embeddings. Any point of accumulation of $(p_t^{(n)})_n$ in $\mathcal{C}([0, 1]^2)$ must equal $p_t^{(0)}$, which follows from testing with suitable functions $\varphi \in L^\infty([0, 1])$. Consequently, $(p_t^{(n)})_n$ is a relatively compact sequence with only one point of accumulation and thus converges. □

PROPOSITION 6.5. *For the class $\Theta_s$ there is a constant $\rho > 0$ such that for all parameters $(\sigma(\cdot), b(\cdot))$ the eigenvalue $\kappa_1 = \kappa_1(\sigma, b)$ of $P_\Delta$ is uniformly separated: $\sigma(P_\Delta) \cap B(\kappa_1, 2\rho) = \{\kappa_1\}$.*

*Furthermore, for all $0 < a < b < 1$ there is a uniform constant $c_{a,b} > 0$ such that the associated first eigenfunction $u_1 = u_1(\sigma, b)$ satisfies, for all $(\sigma, b) \in \Theta_s$,*

$$\min_{x \in [a,b]} |u_1'(x)| \geq c_{a,b}.$$

PROOF. Proceeding indirectly, assume that there is a sequence $(\sigma_n, b_n) \in \Theta_s$ such that the corresponding eigenvalues satisfy $\kappa_1^{(n)} \to 1$ (or $\kappa_1^{(n)} - \kappa_2^{(n)} \to 0$, resp.). By the compactness of the Sobolev embedding of $H^s$ into $C^0$ we can pass to a uniformly converging subsequence. Hence, Proposition 6.4 yields that the corresponding transition densities converge uniformly, which implies that the transition operators $P_\Delta^{(n)}$ converge in operator norm on $L^2([0,1])$. By Proposition 5.6 and Theorem 5.20 in Chatelin (1983), this entails the convergence of their eigenvalues with preservation of the multiplicities. Since the limiting operator is again associated with an elliptic reflected diffusion, the fact that the eigenvalue $\kappa_1 = e^{\Delta \nu_1}$ is always simple (Lemma 6.1) gives the contradiction.

By the same indirect arguments, we construct transition operators $P_\Delta^{(n)}$ on the space $\mathcal{C}([0,1])$ and infer that the eigenfunctions $u_1^{(n)}$ [Chatelin (1983), Theorem 5.10], the invariant measures $\mu^{(n)}$ [see (3.1)] and the inverses of the functions $S^{(n)}$ [see (3.2)] converge in supremum norm. Therefore $(u_1^{(n)})' = \nu_1^{(n)} (S^{(n)})^{-1} \int u_1^{(n)} \mu^{(n)}$ also converges in supremum norm. Due to $u_1'|_{[a,b]} \neq 0$ (Lemma 6.1) this implies that $(u_1^{(n)})'$ cannot converge to zero on $[a,b]$. □

LEMMA 6.6. *The $L^2(\mu)$-normalized eigenfunction $u_k$ of the generator $L$ corresponding to the $(k+1)$st largest eigenvalue $\nu_k$ satisfies*

$$\|u_k\|_{H^{s+1}} \leq C(s, s_0, \|S^{-1}\|_s, \|\mu\|_{s-1}) |\nu_k|^{\lceil s \rceil},$$

*where $C$ is a continuous function of its arguments.*

PROOF. We know that $\mu^{-1}(Su_k')' = \nu_k u_k$ and $u_k'(0) = 0$ holds, which imply

$$u_k'(x) = \nu_k S^{-1}(x) \int_0^x u_k(u) \mu(u)\, du.$$

Suppose $u_k \in H^{r+1}$ with $r \in [0, s]$. Then the function $u_k \mu$ is in $H^{r \wedge (s-1)}$ due to $u_k \in C^r$ (Sobolev embeddings). Hence, the antiderivative is in $H^{(r+1) \wedge s}$. As $S^{-1} \in H^s$ holds, the right-hand side is an element of $H^{s \wedge (r+1)}$. We conclude that the regularity $r$ of $u_k$ is larger by 1, which implies that $u_k$ is in $H^{s+1}$.

In quantitative terms we obtain for $r \in [1, s]$, where we use the seminorm $|f|_s := \|f^{(s)}\|_{L^2}$,

$$\begin{aligned}
|u_k'|_r &\leq |\nu_k| C(r) \left( |S^{-1}|_r \left\| \int_0^\cdot u_k \mu \right\|_\infty + \|S^{-1}\|_\infty \left| \int_0^\cdot u_k \mu \right|_r \right) \\
&\leq |\nu_k| C(r) \|S^{-1}\|_r (\|u_k\|_{L^1(\mu)} + |u_k \mu|_{r-1}) \\
&\leq |\nu_k| C(r) \|S^{-1}\|_s (1 + |u_k|_{r-1} \|\mu\|_\infty + \|u_k\|_\infty \|\mu\|_{r-1}) \\
&\leq |\nu_k| C(r) \|S^{-1}\|_s (1 + 2\|u_k\|_{r-1} \|\mu\|_{s-1}).
\end{aligned}$$

By applying this estimate $\|u_k\|_{r+1} \lesssim |\nu_k| \|S^{-1}\|_s (1 + \|u_k\|_{r-1} \|\mu\|_{s-1})$ successively for $r = 1, 2, \ldots, \lfloor s \rfloor$ and finally, for $r = s - 1$, the estimate follows. $\square$

PROPOSITION 6.7. *For $(\sigma, b) \in \Theta_s$ the corresponding transition probability density $p_\Delta = p_{\Delta, \sigma, b}$ satisfies*

$$\sup_{(\sigma, b) \in \Theta_s} \|p_\Delta\|_{H^{s+1} \times H^s} < \infty.$$

PROOF. The spectral decomposition of $P_\Delta : L^2(\mu) \to L^2(\mu)$ yields

$$p_\Delta(x, y) = \mu(y) \sum_{k=0}^\infty e^{\nu_k \Delta} u_k(x) u_k(y), \qquad x, y \in [0, 1].$$

Due to the uniform ellipticity and uniform boundedness of the coefficients, we have $\nu_k \in [-C_1 k^2, -C_2 k^2]$ with uniform constants $C_1, C_2 > 0$ on $\Theta_s$ [see Davies (1995), adapting Example 4.6.1, page 93, to our situation]. From the preceding Lemma 6.6 and the Sobolev embedding $H^{s+1} \subset C^s$ we infer

$$\begin{aligned}
\|p_\Delta\|_{H^{s+1} \times H^s} &\leq \sum_{k=0}^\infty e^{-C_2 \Delta k^2} \|u_k\|_{s+1} \|\mu u_k\|_s \\
&\leq \sum_{k=0}^\infty C(s, s_0, \|S^{-1}\|_s, \|\mu\|_{s-1})^2 e^{-C_2 \Delta k^2} (C_1 k^2)^{s+1} \|\mu\|_s,
\end{aligned}$$

which gives the desired uniform estimate. $\square$

**Acknowledgment.** The helpful comments of two anonymous referees are gratefully acknowledged.

## REFERENCES

AÏT-SAHALIA, Y. (1996). Nonparametric pricing of interest rate derivative securities. *Econometrica* **64** 527–560.

BANON, G. (1978). Nonparametric identification for diffusion processes. *SIAM J. Control Optim.* **16** 380–395. MR492159

BASS, R. F. (1998). *Diffusions and Elliptic Operators.* Springer, New York. MR1483890

BROWN, B. M. and HEWITT, J. I. (1975). Asymptotic likelihood theory for diffusion processes. *J. Appl. Probability* **12** 228–238. MR375693

CHAPMAN, D. A. and PEARSON, N. D. (2000). Is the short rate drift actually nonlinear? *J. Finance* **55** 355–388.

CHATELIN, F. (1983). *Spectral Approximation of Linear Operators.* Academic Press, New York. MR716134

CHEN, X., HANSEN, L. P. and SCHEINKMAN, J. A. (1997). Shape preserving spectral approximation of diffusions. Working paper. (Last version, November 2000.)

COHEN, A. (2000). Wavelet methods in numerical analysis. In *Handbook of Numerical Analysis* **7** (P. G. Ciarlet, ed.) 417–711. North-Holland, Amsterdam. MR1804747

Davies, E. B. (1995). *Spectral Theory and Differential Operators.* Cambridge Univ. Press. MR1349825

Engel, K.-J. and Nagel, R. (2000). *One-Parameter Semigroups for Linear Evolution Equations.* Springer, Berlin. MR1721989

Fan, J. and Yao, Q. (1998). Efficient estimation of conditional variance functions in stochastic regression. *Biometrika* **85** 645–660. MR1665822

Fan, J. and Zhang, C. (2003). A re-examination of diffusion estimators with applications to financial model validation. *J. Amer. Statist. Assoc.* **98** 118–134. MR1965679

Gobet, E. (2002). LAN property for ergodic diffusions with discrete observations. *Ann. Inst. H. Poincaré Probab. Statist.* **38** 711–737. MR1931584

Hansen, L. P., Scheinkman, J. A. and Touzi, N. (1998). Spectral methods for identifying scalar diffusions. *J. Econometrics* **86** 1–32. MR1650012

Hoffmann, M. (1999). Adaptive estimation in diffusion processes. *Stochastic Process. Appl.* **79** 135–163. MR1670522

Kato, T. (1995). *Perturbation Theory for Linear Operators.* Springer, Berlin. [Reprint of the corrected printing of the 2nd ed. (1980).] MR1335452

Kessler, M. (1997). Estimation of an ergodic diffusion from discrete observations. *Scand. J. Statist.* **24** 211–229. MR1455868

Kessler, M. and Sørensen, M. (1999). Estimating equations based on eigenfunctions for a discretely observed diffusion process. *Bernoulli* **5** 299–314. MR1681700

Kittaneh, F. (1985). On Lipschitz functions of normal operators. *Proc. Amer. Math. Soc.* **94** 416–418. MR787884

Korostelev, A. P. and Tsybakov, A. B. (1993). *Minimax Theory of Image Reconstruction. Lecture Notes in Statist.* **82**. Springer, Berlin. MR1226450

Kutoyants, Y. A. (1975). Local asymptotic normality for processes of diffusion type. *Izv. Akad. Nauk Armyan. SSR Ser. Mat.* **10** 103–112. MR386013

Kutoyants, Y. A. (1984). On nonparametric estimation of trend coefficient in a diffusion process. In *Statistics and Control of Stochastic Processes* 230–250. Nauka, Moscow. MR808204

Liptser, R. S. and Shiryaev, A. N. (2001). *Statistics of Random Processes* **1**. *General Theory*, 2nd ed. Springer, Berlin. MR1800857

Müller, H.-G. and Stadtmüller, U. (1987). Estimation of heteroscedasticity in regression analysis. *Ann. Statist.* **15** 610–625. MR888429

Pham, D. T. (1981). Nonparametric estimation of the drift coefficient in the diffusion equation. *Math. Operationsforsch. Statist. Ser. Statist.* **12** 61–73. MR607955

Reiß, M. (2003). Simulation results for estimating the diffusion coefficient from discrete time observation. Available at www.mathematik.hu-berlin.de/~reiss/sim-diff-est.pdf.

Rosenblatt, M. (1971). *Markov Processes. Structure and Asymptotic Behavior.* Springer, Berlin. MR329037

Stanton, R. (1997). A nonparametric model of term structure dynamics and the market price of interest rate risk. *J. Finance* **52** 1973–2002.

Stroock, D. W. and Varadhan, S. R. S. (1971). Diffusion processes with boundary conditions. *Comm. Pure Appl. Math.* **24** 147–225. MR277037

Stroock, D. W. and Varadhan, S. R. S. (1997). *Multidimensional Diffusion Processes.* Springer, Berlin. MR532498

Tribouley, K. and Viennet, G. (1998). $L_p$ adaptive density estimation in a $\beta$-mixing framework. *Ann. Inst. H. Poincaré Probab. Statist.* **34** 179–208. MR1614587

Yoshida, N. (1992). Estimation for diffusion processes from discrete observations. *J. Multivariate Anal.* **41** 220–242. MR1172898

E. Gobet
Centre de Mathématiques Appliquées
Ecole Polytechnique
91128 Palaiseau Cedex
Paris
France

M. Hoffmann
Laboratoire de Probabilités
   et Modèles Aléatoires
Université Paris 6 et 7
UFR Mathématiques, case 7012 étage
2 Place Jussieu
75251 Paris Cedex 05
France
e-mail: hoffmann@math.jussieu.fr

M. Reiß
Institut für Mathematik
Humboldt-Universität zu Berlin
Unter den Linden 6
10099 Berlin
Germany